\documentclass[english]{article}
\usepackage[latin9]{inputenc}
\usepackage{geometry}
\usepackage{color}
\usepackage{verbatim}
\usepackage{amsmath}
\usepackage{amssymb}
\usepackage{esint}

\makeatletter
\newcommand{\lyxaddress}[1]{
\par {\raggedright #1
\vspace{1.4em}
\noindent\par}
}

\makeatother

\makeatother

\usepackage{babel}

\makeatother

\usepackage{babel}

\makeatother

\usepackage{babel}

\makeatother

\usepackage{babel}

\begin{document}

\title{Integral representations over finite limits for quantum amplitudes}

\author{Jack C. Straton}

\maketitle

\lyxaddress{Department of Physics, Portland State University, Portland, OR 97207-0751,
USA}

\lyxaddress{straton@pdx.edu}
\begin{abstract}
We extend prior work to derive three additional M-1-dimensional integral
representations -- over the interval $[0,1]$, where the prior version
was over the interval $[0,\infty]$ -- for products of M Slater orbitals
(such as appear in quantum transition amplitudes) that allows their
magnitudes of coordinate vector differences (square roots of polynomials)
$|{\bf x}_{1}-{\bf x}_{2}|=\sqrt{x_{1}^{2}-2x_{1}x_{2}\cos\theta+x_{2}^{2}}$
to be moved from disjoint products of functions into a single quadratic
form whose square my be completed. This provides more alternatives
to Fourier transforms that introduce a 3M-dimensional momentum integral
for those products of Slater orbitals (in M separate denominators),
followed in many cases by another set of M-1-dimensional integral
representations to combine those denominators into one denominator
having a single (momentum) quadratic form. The current and prior work
is also slightly more compact than Gaussian transforms that introduce
an M-dimensional integral for products of M Slater orbitals, while
simultaneously moving them into a single (spatial) quadratic form
in a common exponential. 

One may also use addition theorems for extracting the angular variables,
or even direct integration at times. Each method has its strengths
and weaknesses. We have found that two of these M-1-dimensional integral
representations over the interval $[0,1]$ are numerically
stable, as was the prior version having integrals running over the interval
$[0,\infty]$, and one does not need to test for a sufficiently large
upper integration limit as one does for the latter approach. The third integral
representation might have a better form for analytical reduction of integrals. For analytical
reductions of integrals arising from any of the three, however, there is the possible drawback for
large M of there being fewer tabled integrals over $[0,1]$ than over
$[0,\infty]$. In particular, the results of both prior and current
representations have integration variables within square roots as
arguments of Macdonald functions. In a number of cases, these may be
converted to Meijer G-functions whose arguments have the form $(ax^{2}+bx+c)/x$
for which a single tabled integral exists for the integrals from running
over the interval $[0,\infty]$ of the prior paper, and from which
other forms may be found using the techniques given therein. This
is not so for integral representations over the interval $[0,1]$.

Finally, we introduce a fourth integral representation that is not
easily generalizable to large M, but may well provide a bridge for
finding the requisite integrals for such Meijer G-functions over $[0,1]$. 
\end{abstract}
\vspace{2pc}
 \\

\section{Introduction}

When evaluating quantum transition amplitudes, one is faced with the
analytical reduction of integrals involving explicit functions of
the inter-electron (or nucleon) distances. On occasion one may integrate
them directly (see, for instance, \cite{Ley-Koo and Bunge}, among
many others), and at other times addition theorems (e.g., \cite{Sack},
\cite{Porras and  King}, and \cite{Weniger_two-range}) are more
useful. More typically we apply Fourier transforms (e.g., \cite{Fromm and Hill},
\cite{Remiddi}, and \cite{Harris PRA 55 1820}) and/or Gaussian transforms
(e.g., \cite{Kikuchi}, \cite{Shavitt and Karplus}, and \cite{Stra89a})
to effect these reductions.

A prior paper\cite{stra23} introduced a fifth reduction method in
the spirit of Fourier and Gaussian transforms that is an integral
representation having one fewer integral dimension than does a Gaussian
transform  to represent a products of M Slater functions, and roughly
4M fewer integral dimensions than does a Fourier transform  for
such a product. This is an advantage since the main drawback of using
integral representations is that one adds to the number of integral
dimensions one must ultimately solve. In each of these three methods,
the reduction of those introduced integrals becomes more difficult
the larger the numbers of wave functions transformed, so one-fewer
dimension is not a trivial advantage. 

Gaussian transforms require a single one-dimensional integral for
each wave function, and the completion of the square in the coordinate
variables can be done in the resulting exponential. For Fourier transforms,
on the other hand, one must introduce a three-dimensional integral
for each wave function, and often additional integrals to combine
the resulting momentum denominators into a single denominator so that
one can complete the square in the momenta to allow the angular integrals
to be performed.\cite{Stra89c} Our prior work requires the introduction
of one fewer integral dimension than does a Gaussian transform set,
and many fewer than for Fourier transforms. Its main downside is that
the resulting quadratic form (whose square one will complete) resides
in a square root as the argument of a Macdonald function, for which
there are fewer tabled integrals than for the exponential function
wherein resides the quadratic form of Gaussian transforms.

The present paper derives four compact integral representations over
finite intervals to represent a product of Slater-type atomic orbitals,
the seed function $\psi_{000}$ from which Slater functions,\cite{Chen}
Hylleraas powers,\cite{Harris PRA 55 1820} and hydrogenic wave functions
are derived by differentiation. (Known as the Yukawa\cite{Yukawa}
exchange potential in nuclear physics, this function also appears
in plasma physics, where it is known as the Debye-H\"uckel potential,
arising from screened charges\cite{NayekGhoshal} requiring the replacement
of the Coulomb potential by an effective screened potential.\cite{EckerWeizel,Harris}
Such screening of charges also appears in solid-state physics, where
this function is called the Thomas-Fermi potential. In the atomic
physics of negative ions, the radial wave function is given by the
equivalent Macdonald function $\left(R(r)=\frac{C}{\sqrt{r}}K_{1/2}(\eta r)\right)$.\cite{Smirnov2003}
This function also appears in the approximate ground state wave function\cite{GaravelliOliveira}
for a hydrogen atom interacting with hypothesized non-zero-mass photons.\cite{CaccavanoLeung}
We will simply call these \emph{Slater orbitals} herein.)

We start with the simplest integral requiring transformation, the
product of two Slater orbitals integrated over all space,

\begin{equation}
S_{1}^{\eta_{1}0\eta_{12}0}\left(0;0,\mathbf{x}_{2}\right)\equiv S_{1}^{\eta_{1}j_{1}\eta_{12}j_{2}}\left(\mathbf{p}_{1};\mathbf{y}_{1},\mathbf{y}_{2}\right)_{p_{1}\rightarrow0,y_{1},\rightarrow0,y_{2}\rightarrow x_{2},j_{1}\rightarrow0,j_{2}\rightarrow0}=\int d^{3}x_{1}\frac{e^{-\eta_{1}x_{1}}}{x_{1}}\frac{e^{-\eta_{12}x_{12}}}{x_{12}}\quad,\label{eq:SVxVxy}
\end{equation}
where we use the much more general notation of previous work\cite{Stra89a}
in which the short-hand form for shifted coordinates is $\mathbf{x}_{12}=\mathbf{x}_{1}-\mathbf{x}_{2}$,
$\mathbf{p}_{1}$ is a momentum variable within any plane wave associated
with the (first) integration variable, the $\mathbf{y}_{i}$ are coordinates
external to the integration, and the \emph{j}s are defined in the
Gaussian transform\cite{Stra89a} of the generalized Slater orbital:

\begin{equation}
\begin{array}[t]{ccc}
V^{\eta j}({\bf R}) & = & R^{j-1}e^{-\eta R}=\left(-1\right)^{j}{\displaystyle \frac{d^{j}}{d\eta^{j}}}{\displaystyle \frac{1}{\sqrt{\pi}}}\int_{0}^{\infty}\, d\rho{\displaystyle \frac{e^{-R^{2}\rho}e^{-\eta^{2}/(4\rho)}}{\rho^{\;1/2}}}\;\;\left[\eta\geqq0,\: R>0\right]\quad.\\
 & = & \hspace{-0cm}{\displaystyle \frac{1}{2^{j}\sqrt{\pi}}}\int_{0}^{\infty}\, d\rho{\displaystyle \frac{e^{-R^{2}\rho}e^{-\eta^{2}/(4\rho)}}{\rho^{\left(j+1\right)/2}}}H_{j}\left({\displaystyle \frac{\eta}{2\sqrt{\rho}}}\right)\;\;\left[\forall j\geq0\:\mathrm{if\:}\eta>0,\; j=0\:\mathrm{if\:}\eta=0\right]
\end{array}\label{seventeen}
\end{equation}

In our prior work we showed how Gaussian transforms can reduce this
integral in roughly eight steps, so this time we will use Fourier
Transforms:\cite{GR5 p. 512 No. 3.893.1 GR7 p. 486,GR5 p. 382 No. 3.461.2 GR7 p. 364,GR5 p. 384 No. 3.471.9 GR7 p. 368}
\begin{eqnarray}
S_{1}^{\eta_{1}0\eta_{2}0}\left(0;0,x_{2}\right) & = & \int d^{3}x_{1}\frac{e^{-\eta_{1}x_{1}}}{x_{1}}\frac{e^{-\eta_{12}x_{12}}}{x_{12}}\nonumber \\
 & = & \int d^{3}x_{1}\frac{1}{2\pi^{2}}\int\, d^{3}k_{1}\frac{e^{i\,\mathbf{k}_{1}\cdot\mathbf{x}_{\mathrm{1}}}}{\left(\eta_{1}^{2}+k_{1}^{2}\right)}\frac{1}{2\pi^{2}}\int\, d^{3}k_{2}\frac{e^{i\,\mathbf{k}_{2}\cdot{\bf \left(\mathbf{x}_{\mathrm{1}}-\mathbf{x}_{\mathrm{2}}\right)}}}{\left(\eta_{12}^{2}+k_{2}^{2}\right)}\nonumber \\
 & = & \frac{2}{\pi}\int\, d^{3}k_{1}\frac{1}{\left(\eta_{1}^{2}+k_{1}^{2}\right)}\int\, d^{3}k_{2}\frac{e^{-i\,\mathbf{k}_{2}\cdot\mathbf{x}_{\mathrm{2}}}}{\left(\eta_{12}^{2}+k_{2}^{2}\right)}\delta\left(\mathbf{k}_{1}+\mathbf{k}_{2}\right)\nonumber \\
 & = & \frac{2}{\pi}\int\, d^{3}k_{2}\frac{1}{\left(\eta_{1}^{2}+k_{2}^{2}\right)}\frac{e^{-i\,\mathbf{k}_{2}\cdot\mathbf{x}_{\mathrm{2}}}}{\left(\eta_{12}^{2}+k_{2}^{2}\right)}\nonumber \\
 & = & \frac{2}{\pi}\int\, d^{3}k_{2}\int_{0}^{1}\, d\alpha_{1}\frac{1}{\left(\alpha_{1}\left(k_{2}^{2}+\eta_{1}^{2}\right)+\left(1-\alpha_{1}\right)\left(k_{2}^{2}+\eta_{12}^{2}\right)\right)^{2}}e^{-i\,\mathbf{k}_{2}\cdot\mathbf{x}_{\mathrm{2}}}\nonumber \\
 & = & \frac{2}{\pi}\int\, d^{3}k_{2}\int_{0}^{1}\, d\alpha_{1}\int_{0}^{\infty}\, d\rho\rho\exp\left(-i\,\mathbf{k}_{2}\cdot\mathbf{x}_{\mathrm{2}}-\rho\left(\alpha_{1}\left(k_{2}^{2}+\eta_{1}^{2}\right)+\left(1-\alpha_{1}\right)\left(k_{2}^{2}+\eta_{12}^{2}\right)\right)\right)\nonumber \\
 & = & \frac{2}{\pi}\int\, d^{3}k'_{2}\int_{0}^{1}\, d\alpha_{1}\int_{0}^{\infty}\, d\rho\rho\exp\left(-\rho k'{}_{2}^{2}-\frac{x_{2}^{2}}{4\rho}-\rho\left(\alpha_{1}\left(\eta_{1}^{2}-\eta_{12}^{2}\right)+\eta_{12}^{2}\right)\right)\nonumber \\
 & = & \frac{2}{\pi}4\pi\int_{0}^{1}\, d\alpha_{1}\int_{0}^{\infty}\, d\rho\rho\frac{\sqrt{\pi}}{4\rho^{3/2}}\exp\left(-\frac{x_{2}^{2}}{4\rho}-\rho\left(\alpha_{1}\left(\eta_{1}^{2}-\eta_{12}^{2}\right)+\eta_{12}^{2}\right)\right)\nonumber \\
 & = & 2\pi\int_{0}^{1}\, d\alpha_{1}\frac{e^{-x_{2}\sqrt{\alpha_{1}\left(\eta_{1}^{2}-\eta_{12}^{2}\right)+\eta_{12}^{2}}}}{\sqrt{\alpha_{1}\left(\eta_{1}^{2}-\eta_{12}^{2}\right)+\eta_{12}^{2}}}\nonumber \\
 & = & \frac{(4\pi)}{x_{2}\left(\eta_{1}^{2}-\eta_{12}^{2}\right)}\int_{x_{2}\eta_{12}}^{x_{2}\eta_{1}}e^{-y}\, dy\nonumber \\
 & = & \frac{4\pi\left(e^{-\eta_{12}x_{2}}-e^{-\eta_{1}x_{2}}\right)}{x_{2}\left(\eta_{1}^{2}-\eta_{12}^{2}\right)}\quad.\label{eq:YYFourier}
\end{eqnarray}

This is a considerably lengthy derivation, and the Gaussian Transform
approach is not much better. Of course, in this simple case, one can
invoke the addition theorem expression for $e^{-i\,\mathbf{k}_{2}\cdot\mathbf{x}_{\mathrm{2}}}$
in the fourth line to shorten the reduction, but for the large-M equivalent,
the above process is what one must follow . (Actually, one would do
well to avoid the use of the Dirac delta function in the third line
when dealing with large M.)

This was the motivation for the fifth path to a solution of our prior
paper.

\section{A simpler integral representation}

We begin by introducing an integral representation over a finite interval  for a pair-product
of Slater orbitals,

\begin{align}
\frac{e^{-\eta_{1}x_{1}}}{x_{1}}\frac{e^{-\eta_{12}x_{12}}}{x_{12}} & =\int_{0}^{1}\, d\alpha_{1}\frac{\sqrt{\left(1-\alpha_{1}\right)\eta_{1}^{2}+\alpha_{1}\eta_{12}^{2}}K_{1}\left(\sqrt{\frac{x_{1}^{2}}{1-\alpha_{1}}+\frac{x_{12}^{2}}{\alpha_{1}}}\sqrt{\left(1-\alpha_{1}\right)\eta_{1}^{2}+\alpha_{1}\eta_{12}^{2}}\right)}{\pi\left(1-\alpha_{1}\right)^{3/2}\alpha_{1}^{3/2}\sqrt{\frac{x_{1}^{2}}{1-\alpha_{1}}+\frac{x_{12}^{2}}{\alpha_{1}}}}\quad,\label{eq:YYalpha}
\end{align}
whose derivation will follow a display of its utility. We insert it
in the above problem, %
{} completing the square in the quadratic form (changing variables from
$\mathbf{x}_{1}$ to $\mathbf{x}'_{1}=\mathbf{x}_{1}-\left(1-\alpha_{1}\right)\mathbf{x}_{2}$
with unit Jacobian) in both places in the integrand where it appears,
so that \cite{GR5 p. 727 No. 6.596.3 GR7 p. 693,GR5 p. 111 No. 2.311}

\begin{eqnarray}
S_{1}^{\eta_{1}0\eta_{12}0}\left(0;0,x_{2}\right) & = & \int d^{3}x_{1}\int_{0}^{1}\, d\alpha_{1}\frac{\sqrt{\left(1-\alpha_{1}\right)\eta_{1}^{2}+\alpha_{1}\eta_{12}^{2}}K_{1}\left(\sqrt{\frac{x_{1}^{2}}{1-\alpha_{1}}+\frac{x_{12}^{2}}{\alpha_{1}}}\sqrt{\left(1-\alpha_{1}\right)\eta_{1}^{2}+\alpha_{1}\eta_{12}^{2}}\right)}{\pi\left(1-\alpha_{1}\right)^{3/2}\alpha_{1}^{3/2}\sqrt{\frac{x_{1}^{2}}{1-\alpha_{1}}+\frac{x_{12}^{2}}{\alpha_{1}}}}\nonumber \\
 & = & \int d^{3}x'_{1}\int_{0}^{1}\, d\alpha_{1}\frac{\sqrt{\left(1-\alpha_{1}\right)\eta_{1}^{2}+\alpha_{1}\eta_{12}^{2}}K_{1}\left(\sqrt{\frac{x_{1}^{'2}}{\left(1-\alpha_{1}\right)\alpha_{1}}+x_{2}^{2}}\sqrt{\left(1-\alpha_{1}\right)\eta_{1}^{2}+\alpha_{1}\eta_{12}^{2}}\right)}{\pi\left(1-\alpha_{1}\right)^{3/2}\alpha_{1}^{3/2}\sqrt{\frac{x_{1}^{'2}}{\left(1-\alpha_{1}\right)\alpha_{1}}+x_{2}^{2}}}\nonumber \\
 & = & \int_{0}^{1}\, d\alpha_{1}\frac{2\pi e^{-x_{2}\sqrt{\left(1-\alpha_{1}\right)\eta_{1}^{2}+\alpha_{1}\eta_{12}^{2}}}}{\sqrt{\left(1-\alpha_{1}\right)\eta_{1}^{2}+\alpha_{1}\eta_{12}^{2}}}=\int_{x_{2}\eta_{1}}^{x_{2}\eta_{12}}\, dy\frac{4\pi e^{-y}}{x_{2}\left(\eta_{12}^{2}-\eta_{1}^{2}\right)}\nonumber \\
 & = & \frac{4\pi\left(e^{-\eta_{12}x_{2}}-e^{-\eta_{1}x_{2}}\right)}{x_{2}\left(\eta_{1}^{2}-\eta_{12}^{2}\right)}\quad,\label{eq:syyvianew-1}
\end{eqnarray}
which is indeed a much shorter path to the solution than the Fourier
and Gaussian transforms give. This follows from the fact that it requires
the introduction of one integral to represent a pair-product of Slater
orbitals rather than one integral for each orbital that the Gaussian
transform requires, or the three-dimensional integral that the Fourier
transform approach requires for each Slater orbital (with two additional
integrals required, as in eq. (\ref{eq:YYFourier}) ).

Note that this new integral representation has a similar integrand
to the integral representation introduced in our prior paper for 
products of M Slater orbitals,

\begin{equation}
\begin{array}{ccc}
\frac{e^{-R_{1}\eta_{1}}}{R_{1}}\hspace{-0.3cm} & \cdot & \hspace{-0.6cm}\frac{e^{-R_{2}\eta_{2}}}{R_{2}}\cdots\frac{e^{-R_{M}\eta_{M}}}{R_{M}}=\frac{1}{2^{M}\pi^{2M}}\int_{0}^{\infty}d\zeta_{1}\int_{0}^{\infty}d\zeta_{2}\cdots\int_{0}^{\infty}\, d\zeta_{M-1}{\displaystyle \frac{\pi^{3M/2}}{\prod_{i=1}^{M-1}\zeta_{i}^{3/2}}}2^{\frac{M}{2}+1}\\
 & \times & \left(R_{1}^{2}+\frac{R_{2}^{2}}{\zeta_{1}}+\frac{R_{3}^{2}}{\zeta_{2}}+\cdots\,+\frac{R_{M}^{2}}{\zeta_{M-1}}\right)^{-M/4}\left(\eta_{1}^{2}+\zeta_{1}\eta_{2}^{2}+\zeta_{2}\eta_{3}^{2}+\cdots\,+\zeta_{M-1}\eta_{M}^{2}\right)^{M/4}\\
 & \times & K_{\frac{M}{2}}\left(\sqrt{R_{1}^{2}+\frac{R_{2}^{2}}{\zeta_{1}}+\frac{R_{3}^{2}}{\zeta_{2}}+\cdots\,+\frac{R_{M}^{2}}{\zeta_{M-1}}}\sqrt{\eta_{1}^{2}+\zeta_{1}\eta_{2}^{2}+\zeta_{2}\eta_{3}^{2}+\cdots\,+\zeta_{M-1}\eta_{M}^{2}}\right)\quad,
\end{array}\label{eq:MtransCompact}
\end{equation}
except that the new integral representation has finite limits of integration
rather than the infinite interval of the previous work.

\section{Problematic Approaches}

The first step in creating the prior integral representation -- and
the new one -- entails converting a product of Slater orbitals into
denominators of some power (combined with other factors) using some
\emph{initial} integral representation, such as via the Stieltjes
Transform \cite{ET II p. 220 No. 14.2.41}, or \cite{GR5 p. 706 No. 6.554.4 GR7 p. 675},
or \cite{GR5 p. 467 No. 3.773.5 GR7 p. 444}, below

\begin{equation}
\frac{e^{-\eta_{1}x_{1}}}{x_{1}}\frac{e^{-\eta_{12}x_{12}}}{x_{12}}=\int_{0}^{\infty}\int_{0}^{\infty} dt_{1}dt_{2} \frac{2}{\pi}\frac{\cos\left(t_{1}\eta_{1}\right)}{\left(t_{1}^{2}+x_{1}^{2}\right)}\frac{2}{\pi}\frac{\cos\left(t_{2}\eta_{12}\right)}{\left(t_{2}^{2}+x_{12}^{2}\right)} \quad. \label{eq:zt1t2_trans}
\end{equation}
In the prior paper we combined products of denominators into one,
consolidating the coordinate variables into a common quadratic form,
using \cite{GR5 p. 649 No. 4.638.2 GR7 p. 615}

\begin{eqnarray}
\frac{1}{r_{1}^{p_{1}}r_{2}^{p_{2}}\cdots r_{n}^{p_{n}}r_{0}^{s-p_{1}-p_{2}-\cdots-p_{n}}} & = & \frac{\Gamma\left(s\right)}{\Gamma\left(s-p_{1}-p_{2}-\cdots-p_{n}\right)\Gamma\left(p_{1}\right)\Gamma\left(p_{2}\right)\cdots\Gamma\left(p_{n}\right)}
\nonumber \\
& \times & \int_{0}^{\infty}d\zeta_{1}\int_{0}^{\infty}d\zeta_{2}\cdots\int_{0}^{\infty}d\zeta_{n}\frac{\zeta_{1}^{p_{1}-1}\zeta_{1}^{p_{2}-1}\cdots\zeta_{1}^{p_{n}-1}}{\left(r_{0}+r_{1}\zeta_{1}+r_{2}\zeta_{2}+\cdots+r_{n}\zeta_{n}\right)^{s}}\nonumber \\
 & & \quad\left[p_{i}>0,\, r_{i}>0,\, s>p_{1}+p_{2}+\cdots+p_{n}>0\right]\quad,\label{eq:GR5 p. 649 No. 4.638.2 GR7 p. 615}
\end{eqnarray}
(with $p_{1}=1$ and $s=2$ for $n=1$):

\begin{equation}
\frac{e^{-\eta_{1}x_{1}}}{x_{1}}\frac{e^{-\eta_{12}x_{12}}}{x_{12}}=\frac{4}{\pi^{2}}\int_{0}^{\infty}\int_{0}^{\infty}dt_{1}dt_{2} \cos\left(t_{1}\eta_{1}\right)\cos\left(t_{2}\eta_{12}\right)\int_{0}^{\infty} d\zeta_{1} \frac{\zeta_{1}^{1-1}}{\left(\zeta_{1}\left(t_{1}^{2}+x_{1}^{2}\right)+\left(t_{2}^{2}+x_{12}^{2}\right)\right)^{2}}\quad.\label{eq:t1t2zeta-1}
\end{equation}
One then performs the\emph{ t} integrals to obtain the $n=2$ version
of eq. (\ref{eq:MtransCompact}).

To represent a product of M Slater orbitals using finite-interval
integrals, one can in principle use Feynman parametrization\cite{Feynman 1949}  as extended
by Schweber.\cite{Schweber} His third version of the extension, which
may be derived by iterating \cite{GH p. 175 No. 421.4 ,GR5 p. 336 No. 3.199 GR7 p. 318} is 

\begin{eqnarray}
\hspace{-0.5cm}   \frac{1}{D_{1}D_{2}\cdots D_{n}}  \hspace{-0.3cm} & = & \hspace{-0.3cm}  \left(n-1\right)!\int_{0}^{1}d\alpha_{1}\alpha_{1}^{n-2}\int_{0}^{1}d\alpha_{2}\alpha_{2}^{n-3}\cdots\int_{0}^{1}d\alpha_{n-1}
\nonumber \\
& \times & \hspace{-0.3cm}  \frac{1}{\left(D_{1}\alpha_{1}\alpha_{2}\cdots\alpha_{n-1}+D_{2}\alpha_{1}\cdots\alpha_{n-2}\left(1-\alpha_{n-1}\right)+\cdots+D_{n-1}\alpha_{1}\left(1-\alpha_{2}\right)+D_{n}\left(1-\alpha_{1}\right)\right)^{n}}
\end{eqnarray}
so that, for instance,

\begin{equation}
\frac{e^{-\eta_{1}x_{1}}}{x_{1}}\frac{e^{-\eta_{12}x_{12}}}{x_{12}}=\frac{4}{\pi^{2}}\int_{0}^{\infty}\int_{0}^{\infty}\cos\left(t_{1}\eta_{1}\right)\cos\left(t_{2}\eta_{12}\right) dt_{1}dt_{2} \int_{0}^{1} d\alpha_{1} \frac{1}{\left(\left(t_{1}^{2}+x_{1}^{2}\right)\alpha_{1}+\left(t_{2}^{2}+x_{12}^{2}\right)\left(1-\alpha_{1}\right)\right)^{2}} \, ,\label{eq:t1t2zeta-1-1}
\end{equation}
and this indeed bears fruit and is numerically stable up through $M=3$.
But the subsequent integration of the \emph{t}s (as in the $n=2$
version eq. (\ref{eq:zt1t2_trans})) produces functions whose arguments
do not form discernable patterns as \emph{M} increases, so we were
not able to generalize this approach (using  \cite{GR5 p. 467 No. 3.773.5 GR7 p. 444}) to large
\emph{M}.

In our prior work we used Fourier transforms to convert the product
of M Slater orbitals into denominators instead of the integral set
using \cite{GR5 p. 467 No. 3.773.5 GR7 p. 444} that we utilized in
eq. (\ref{eq:zt1t2_trans}). On a formal level, this approach worked
fine with Schweber's third parametrization, yielding

\begin{eqnarray}
\frac{e^{-R_{1}\eta_{1}}}{R_{1}}\cdot\frac{e^{-R_{2}\eta_{2}}}{R_{2}}\cdots\frac{e^{-R_{M}\eta_{M}}}{R_{M}} & = & \frac{1}{2^{M}\pi^{2M}}\int_{0}^{\infty}d\rho\int_{0}^{1}d\alpha_{1}\int_{0}^{1}d\alpha_{2}\cdots\int_{0}^{1}\, d\alpha_{M-1}\alpha_{1}^{M-2}\alpha_{2}^{M-3}\cdots\alpha_{M-3}^{2}\alpha_{M-2}^{1}
\nonumber \\
& \times & {\displaystyle \frac{\pi^{3M/2}}{\rho^{M/2+1}\prod_{i=1}^{M-1}\alpha_{i}^{3\left(M-i\right)/2}\left(1-\alpha_{i}\right)^{3/2}}}\\
 & \times & exp\left(-\rho\left(\alpha_{1}\alpha_{2}\cdots\alpha_{M-1}\eta_{1}^{2}+\alpha_{1}\alpha_{2}\cdots\alpha_{M-2}\left(1-\alpha_{M-1}\right)\eta_{2}^{2} \right.  \right.
 \nonumber \\  & + & \left.  \left. \cdots+\alpha_{1}\left(1-\alpha_{2}\right)\eta_{M-1}^{2}+\cdots\,+\left(1-\alpha_{1}\right)\eta_{M}^{2}\right)\right)\\
 & \times & exp\left(- \left(\frac{R_{1}^{2}}{\alpha_{1}\alpha_{2}\cdots\alpha_{M-1}}+\frac{R_{2}^{2}}{\alpha_{1}\alpha_{2}\cdots\alpha_{M-2}\left(1-\alpha_{M-1}\right)}\right.  \right.
 \nonumber \\
  & + & \left.  \left. \cdots+\frac{R_{M-1}^{2}}{\alpha_{1}\left(1-\alpha_{2}\right)}\,+\frac{R_{M}^{2}}{\left(1-\alpha_{1}\right)}\right)\frac{1}{4\rho}\right)\quad.
\label{eq:MRpoForm-1-1}
\end{eqnarray}
However, this was numerically stable only up through $M=3$. After
performing the $\rho$ integral \cite{GR5 p. 384 No. 3.471.9 GR7 p. 368} to
give the most compact form, the result was numerically stable only
up through $M=2$:

\begin{equation}
\begin{array}{ccc}
& & \hspace{-1.9cm}  \frac{e^{-R_{1}\eta_{1}}}{R_{1}} \cdot  \frac{e^{-R_{2}\eta_{2}}}{R_{2}}\cdots\frac{e^{-R_{M}\eta_{M}}}{R_{M}}=\frac{1}{2^{M}\pi^{2M}}\int_{0}^{1}d\alpha_{1}\int_{0}^{1}d\alpha_{2}\cdots\int_{0}^{1}\, d\alpha_{M-1}{\displaystyle \alpha_{1}^{M-2}\alpha_{2}^{M-3}\cdots\alpha_{M-3}^{2}\alpha_{M-2}^{1}\frac{\pi^{3M/2}}{\prod_{i=1}^{M-1}\zeta_{i}^{3/2}}}2^{\frac{M}{2}+1}\\
 & \times &  \hspace{-1.1cm} \left(\alpha_{1}\alpha_{2}\cdots\alpha_{M-1}\eta_{1}^{2}+\alpha_{1}\alpha_{2}\cdots\alpha_{M-2}\left(1-\alpha_{M-1}\right)\eta_{2}^{2}+\cdots+\alpha_{1}\left(1-\alpha_{2}\right)\eta_{M-1}^{2}+\cdots\,+\left(1-\alpha_{1}\right)\eta_{M}^{2}\right)^{M/4}\\
 & \times &  \hspace{-5.9cm} \left(\frac{R_{1}^{2}}{\alpha_{1}\alpha_{2}\cdots\alpha_{M-1}}+\frac{R_{2}^{2}}{\alpha_{1}\alpha_{2}\cdots\alpha_{M-2}\left(1-\alpha_{M-1}\right)}+\cdots+\frac{R_{M-1}^{2}}{\alpha_{1}\left(1-\alpha_{2}\right)}\,+\frac{R_{M}^{2}}{\left(1-\alpha_{1}\right)}\right)^{-M/4}\\
 & \times & K_{\frac{M}{2}}\left(\sqrt{\alpha_{1}\alpha_{2}\cdots\alpha_{M-1}\eta_{1}^{2}+\alpha_{1}\alpha_{2}\cdots\alpha_{M-2}\left(1-\alpha_{M-1}\right)\eta_{2}^{2}+\cdots+\alpha_{1}\left(1-\alpha_{2}\right)\eta_{M-1}^{2}+\cdots\,+\left(1-\alpha_{1}\right)\eta_{M}^{2}}\right.\quad \\
 & \times & \left.\sqrt{\frac{R_{1}^{2}}{\alpha_{1}\alpha_{2}\cdots\alpha_{M-1}}+\frac{R_{2}^{2}}{\alpha_{1}\alpha_{2}\cdots\alpha_{M-2}\left(1-\alpha_{M-1}\right)}+\cdots+\frac{R_{M-1}^{2}}{\alpha_{1}\left(1-\alpha_{2}\right)}\,+\frac{R_{M}^{2}}{\left(1-\alpha_{1}\right)}}\right) \,.
\end{array}\label{eq:MtransCompact-Schweber3}
\end{equation}

This approach, then, should be reserved for analytical reduction of
integrals rather than numerical integration.

\section{Derivation of a Second, Numerically-Stable Integral Representation}

Schweber's second parametrization,

\begin{eqnarray}
\hspace{-0.5cm}   \frac{1}{D_{1}D_{2}\cdots D_{n}} & = & \left(n-1\right)!\int_{0}^{1}d\alpha_{1}\int_{0}^{\alpha_{1}}d\alpha_{2}\cdots\int_{0}^{\alpha_{n-2}}d\alpha_{n-1}
\nonumber \\
& \times & \hspace{-0.3cm} 
\frac{1}{\left(D_{n}\alpha_{n-1}+D_{n-1}\left(\alpha_{n-2}-\alpha_{n-1}\right)+\cdots+D_{1}\left(1-\alpha_{1}\right)\right)^{n}}
\end{eqnarray}

\noindent
looks somewhat dubious for analytical uses since each succeeding integral
has the prior parameter as its upper limit. It turns out, however,
to give numerically stable results (we checked up through $M=6$ when using the Fourier transform as the bridge
rather than (\ref{eq:zt1t2_trans})), and one may perform a change
of variables in each integral at the end of the derivation to give all
integrals over {[}0,1{]}. Since Fourier transforms include momentum
variables in plane waves, we take the additional step of moving the
combined momentum denominator into an exponential by using \cite{GR5 p. 364 No. 3.381.4}
\begin{equation}
(\nu-1)!D^{-\nu}=\int_{0}^{\infty}d\rho\rho^{\nu-1}e^{-\rho D}\quad.\label{eq:8}
\end{equation}

Thus, for a product of M Slater orbitals, we have\cite{GR5 p. 649 No. 4.638.2 GR7 p. 615}

\begin{eqnarray}
\hspace{-0.25cm}  \label{eq:Schweber3FourierApproach} \frac{e^{-R_{1}\eta_{1}}}{R_{1}}\hspace{-0.45cm} & \cdot & \hspace{-0.45cm}\frac{e^{-R_{2}\eta_{2}}}{R_{2}}\cdots\frac{e^{-R_{M}\eta_{M}}}{R_{M}}=\int d^{3}k_{1}\int d^{3}k_{2}\cdots\int\, d^{3}k_{M}\frac{1}{2\pi^{2}}\cdot\frac{e^{ik_{1}\cdot R_{1}}}{k_{1}^{2}+\eta_{1}^{2}}\cdot\frac{1}{2\pi^{2}}\cdot\frac{e^{ik_{2}\cdot R_{2}}}{k_{2}^{2}+\eta_{2}^{2}}\cdots\,\frac{1}{2\pi^{2}}\cdot\frac{e^{ik_{M}\cdot R_{M}}}{k_{M}^{2}+\eta_{M}^{2}}\nonumber \\
 & = &  \hspace{-0.35cm}   \int_{0}^{1}d\alpha_{1}\int_{0}^{\alpha_{1}}d\alpha_{2}\cdots\int_{0}^{\alpha_{M-2}}\, d\alpha_{M-1}\int d^{3}k_{1}\int d^{3}k_{2}\cdots\int\, d^{3}k_{M}\frac{(M-1)!}{2^{M}\pi^{2M}} \nonumber \\
 & \times & \hspace{-0.45cm} \frac{\exp\left(ik_{1}\cdot R_{1}+ik_{2}\cdot R_{2}+\cdots+ik_{M-1}\cdot x_{M-1}+ik_{M}\cdot R_{M}\right)}{\left(\left(k_{1}^{2}+\eta_{1}^{2}\right)\left(1-\alpha_{1}\right)+\left(k_{2}^{2}+\eta_{2}^{2}\right)\left(\alpha_{1}-\alpha_{2}\right)  \hspace{-0.10cm} + \hspace{-0.10cm} \cdots \hspace{-0.08cm} + \hspace{-0.10cm} \left(k_{M-1}^{2}+\eta_{M-1}^{2}\right)\left(\alpha_{M-2}-\alpha_{M-1}\right)+\left(k_{M}^{2}+\eta_{M}^{2}\right)\alpha_{M}\right)^{M}}  \nonumber \\
 & = &  \hspace{-0.35cm}  \frac{1}{2^{M}\pi^{2M}}\int_{0}^{\infty}d\rho\int_{0}^{1}d\alpha_{1}\int_{0}^{\alpha_{1}}d\alpha_{2}\cdots\int_{0}^{\alpha_{M-2}}\, d\alpha_{M-1}\int d^{3}k_{1}\int d^{3}k_{2}\cdots\int\, d^{3}k_{M}  \rho^{M-1} \nonumber \\
 & \times &  \hspace{-0.35cm}  \exp\left(-\rho\left(ik_{1}\cdot R_{1}/\rho+ik_{2}\cdot R_{2}/\rho+\cdots+ik_{M-1}\cdot x_{M-1}/\rho+ik_{M}\cdot R_{M}/\rho\right)\right)
 exp\left(-\rho\left(\left(k_{1}^{2}+\eta_{1}^{2}\right)\left(1-\alpha_{1}\right) \right. \right. \nonumber \\
 & +&  \hspace{-0.35cm}  \left. \left. \left(k_{2}^{2}+\eta_{2}^{2}\right)\left(\alpha_{1}-\alpha_{2}\right)+\cdots+\left(k_{M-1}^{2}+\eta_{M-1}^{2}\right)\left(\alpha_{M-2}-\alpha_{M-1}\right)+\left(k_{M}^{2}+\eta_{M}^{2}\right)\alpha_{M}\right)\right)\nonumber \\
 & \equiv &  \hspace{-0.35cm}  \frac{1}{2^{M}\pi^{2M}}\int_{0}^{\infty}d\rho\int_{0}^{1}d\alpha_{1}\int_{0}^{\alpha_{1}}d\alpha_{2}\cdots\int_{0}^{\alpha_{M-2}}\, d\alpha_{M-1}\int d^{3}k_{1}\int d^{3}k_{2}\cdots\int\, d^{3}k_{M}\rho^{M-1}\exp\left(-\rho\, Q\right)\,. \label{eq:Schweber3FourierApproach} 
\end{eqnarray}

The quadratic form may be written as\cite{Stra89c}

\begin{equation}
Q=\underline{V}^{T}\underline{W}\underline{V}\quad,\label{eq:10}
\end{equation}
 where 
\begin{equation}
\underline{V}^{T}=\left(\mathbf{k}_{1},\,\mathbf{k}_{2},\cdots\mathbf{,\, k}_{M},1\right)\quad,\label{eq:11}
\end{equation}

\begin{equation}
\underline{W}=\left(\begin{array}{cccccc}
\left(1-\alpha_{1}\right) & 0 & \cdots & 0 & 0 & {\bf b}_{1}\\
0 & \left(\alpha_{1}-\alpha_{2}\right) & \cdots & 0 & 0 & {\bf b}_{2}\\
\vdots & \vdots & \ddots & \vdots & \vdots & \vdots\\
0 & 0 & \cdots & \left(\alpha_{M-2}-\alpha_{M-1}\right) & 0 & {\bf b}_{M-1}\\
0 & 0 & \cdots & 0 & \alpha_{M-1} & {\bf b}_{M}\\
{\bf b}_{1} & {\bf b}_{2} & \cdots & {\bf b}_{M-1} & {\bf b}_{M} & C
\end{array}\right)\;\;,\label{eq:12}
\end{equation}

\begin{equation}
C=\left(1-\alpha_{1}\right)\eta_{1}^{2}+\left(\alpha_{1}-\alpha_{2}\right)\eta_{2}^{2}+\cdots+\left(\alpha_{M-2}-\alpha_{M-1}\right)\eta_{M-1}^{2}+\alpha_{M-1}\eta_{M}^{2}\;\;,\label{eq:13}
\end{equation}
 and

\begin{equation}
\mathbf{b_{j}}=-\frac{i}{2\rho}\mathbf{R_{j}}\quad.\label{eq:15}
\end{equation}

Now suppose one could find an orthogonal transformation that reduced
Q to diagonal form

\begin{equation}
Q'=a'_{1}k{}_{1}^{'2}+a'_{'2}k_{2}^{'2}+\ldots+a'_{N+M}k_{N+M}^{'2}+c',\label{eq:16}
\end{equation}
where, as shown by Chisholm,\cite{Chisholm} the $a'$ are positive.
Then after a simple translation in $\lbrace\mathbf{k}_{1},\,\mathbf{k}_{2},\cdots,\mathbf{\, k}_{M}\rbrace$
space (with Jacobian = 1), the \emph{k} integrals could be done,\cite{GR5 p. 382 No. 3.461.2 GR7 p. 364}
\begin{equation}
\int\, d^{3}k'_{1}\ldots\, d^{3}k'_{M}e^{-\rho\left(a'_{1}k{}_{1}^{'2}+a'_{2}k_{2}^{'2}+\ldots+a'_{M}k_{M}^{'2}\right)}=\left(\frac{\pi^{M}}{\rho^{M}\Lambda}\right)^{3/2}.\label{eq:17}
\end{equation}
 Since this result is expressed in the form of an invariant determinant,

\begin{equation}
\Lambda=
\hspace{-0.10cm}
\left|\begin{array}{ccccc}
\hspace{-0.10cm} \left(1-\alpha_{1}\right) & 0 & \cdots & 0 & 0\\
0 & \left(\alpha_{1}-\alpha_{2}\right) & \cdots & 0 & 0\\
\vdots & \vdots & \ddots & \vdots & \vdots\\
0 & 0 & \cdots & \left(\alpha_{M-2}-\alpha_{M-1}\right) & 0\\
\hspace{-0.10cm}0 & 0 & \cdots & 0 & \alpha_{M-1}
\end{array} \hspace{-0.10cm} \right|=\left(1-\alpha_{1}\right)\left(\prod_{i=2}^{M-1}\left(\alpha_{i-1}-\alpha_{i}\right)\right)\alpha_{M-1}=\prod_{i=1}^{M}a'_{i}\;,\label{eq:18}
\end{equation}
actually finding the orthogonal transformation that reduces Q to diagonal
form is unnecessary. What is left to find is just the exponential
of $-\rho c'$, which we integrate over $\rho$ and the $\alpha_{i}$.

This orthogonal transformation also leaves

\begin{equation}
\Omega=\mathrm{det}\mathbf{W}\label{eq:35}
\end{equation}

\noindent invariant and to find its value one need only expand $\Omega$
by minors:

\begin{eqnarray}
\Omega & = & C\Lambda+\sum_{i=1}^{M}\sum_{j=1}^{M}\mathrm{\mathbf{b}}_{i}\cdot\mathrm{\mathbf{b}}_{j}\left(-1\right)^{i+j+1}\Lambda_{ij} \nonumber  \\
 & = & C\Lambda-b_{1}^{2}\left(\prod_{i=2}^{M-1}\left(\alpha_{i-1}-\alpha_{i}\right)\right)\alpha_{M-1}-\sum_{j=2}^{M-1}b_{j}^{2}\frac{\left(1-\alpha_{1}\right)}{\alpha_{j-1}-\alpha_{j}}\left(\prod_{i=2}^{M-1}\left(\alpha_{i-1}-\alpha_{i}\right)\right)\alpha_{M-1} \nonumber \\
 &- &b_{M}^{2}\left(1-\alpha_{1}\right)\prod_{i=2}^{M-1}\left(\alpha_{i-1}-\alpha_{i}\right)\quad,\label{eq:40}
\end{eqnarray}

\noindent where $\Lambda_{ij}$ is $\Lambda$ with the \emph{i}th
row and \emph{j}th column deleted, and  is diagonal in the present
case. Therefore, $c'$ (of eq. (\ref{eq:16})) is given by

\begin{eqnarray}
c'=\Omega/\Lambda & = & \left(1-\alpha_{1}\right)\eta_{1}^{2}+\sum_{j=2}^{M-1}\eta_{j}^{2}\left(\alpha_{j-1}-\alpha_{j}\right)+\alpha_{M-1}\eta_{M}^{2}-\frac{b_{1}^{2}}{\left(1-\alpha_{1}\right)}-\sum_{j=2}^{M-1}\frac{b_{j}^{2}}{\left(\alpha_{j-1}-\alpha_{j}\right)}-\frac{b_{M}^{2}}{\alpha_{M-1}}\nonumber \\
 & = & \left(1-\alpha_{1}\right)\eta_{1}^{2}+\sum_{j=2}^{M-1}\eta_{j}^{2}\left(\alpha_{j-1}-\alpha_{j}\right)+\alpha_{M-1}\eta_{M}^{2}+\frac{1}{4\rho^{2}}\frac{R_{1}^{2}}{\left(1-\alpha_{1}\right)}+\sum_{j=2}^{M-1}\frac{1}{4\rho^{2}}\frac{R_{j}^{2}}{\left(\alpha_{j-1}-\alpha_{j}\right)} \nonumber \\
 &+& \frac{1}{4\rho^{2}}\frac{b_{M}^{2}}{\alpha_{M-1}} \quad , \label{eq:39}
\end{eqnarray}
 so that

\begin{equation}
\begin{array}{ccc}
\frac{e^{-R_{1}\eta_{1}}}{R_{1}} & \cdot &  \hspace{-2.6cm}  \frac{e^{-R_{2}\eta_{2}}}{R_{2}}\cdots\frac{e^{-R_{M}\eta_{M}}}{R_{M}}  =  \frac{1}{2^{M}\pi^{2M}}\int_{0}^{\infty}d\rho\int_{0}^{1}d\alpha_{1}\int_{0}^{\alpha_{1}}d\alpha_{2}\cdots\int_{0}^{\alpha_{M-2}}\, d\alpha_{M-1}\\
 & \times &   \hspace{-5.5cm}  {\displaystyle \frac{\pi^{3M/2}}{\rho^{M/2+1}\left(\left(1-\alpha_{1}\right)\left(\prod_{i=2}^{M-1}\left(\alpha_{i-1}-\alpha_{i}\right)\right)\alpha_{M-1}\right)^{3/2}}}\\
 & \times & exp\left(-\rho\left(\left(1-\alpha_{1}\right)\eta_{1}^{2}+\left(\alpha_{1}-\alpha_{2}\right)\eta_{2}^{2}+\cdots+\left(\alpha_{M-2}-\alpha_{M-1}\right)\eta_{M-1}^{2}+\cdots\,+\alpha_{M-1}\eta_{M}^{2}\right)\right)\\
 & \times &  \hspace{-1.5cm} exp\left(-{\displaystyle \left(\frac{R_{1}^{2}}{\left(1-\alpha_{1}\right)}+\frac{R_{2}^{2}}{\left(\alpha_{1}-\alpha_{2}\right)}+\cdots+\frac{R_{M-1}^{2}}{\left(\alpha_{M-2}-\alpha_{M-1}\right)}\,+\frac{R_{M}^{2}}{\alpha_{M-1}}\right)\frac{1}{4\rho}}\right)\quad.
\end{array}\label{eq:MRpoForm-Schweber2}
\end{equation}
 We perform the $\rho$ integral \cite{GR5 p. 384 No. 3.471.9 GR7 p. 368} to
give the most compact, semi-final form for the desired integral representation:

\begin{equation}
\begin{array}{ccc}
\frac{e^{-R_{1}\eta_{1}}}{R_{1}}\hspace{-0.3cm} & \cdot & \hspace{-0.6cm}\frac{e^{-R_{2}\eta_{2}}}{R_{2}}\cdots\frac{e^{-R_{M}\eta_{M}}}{R_{M}}=\int_{0}^{1}d\alpha_{1}\int_{0}^{\alpha_{1}}d\alpha_{2}\cdots\int_{0}^{\alpha_{M-2}}\, d\alpha_{M-1}{\displaystyle \frac{2^{1-\frac{M}{2}}\pi^{-M/2}}{\left(1-\alpha_{1}\right)^{3/2}\left(\prod_{i=2}^{M-1}\left(\alpha_{i-1}-\alpha_{i}\right)^{3/2}\right)\alpha_{M-1}^{3/2}}}\\
 & \times & \left(\left(1-\alpha_{1}\right)\eta_{1}^{2}+\left(\alpha_{1}-\alpha_{2}\right)\eta_{2}^{2}+\cdots+\left(\alpha_{M-2}-\alpha_{M-1}\right)\eta_{M-1}^{2}+\cdots\,+\alpha_{M-1}\eta_{M}^{2}\right)^{M/4}\\
 & \times & \left(\frac{R_{1}^{2}}{\left(1-\alpha_{1}\right)}+\frac{R_{2}^{2}}{\left(\alpha_{1}-\alpha_{2}\right)}+\cdots+\frac{R_{M-1}^{2}}{\left(\alpha_{M-2}-\alpha_{M-1}\right)}\,+\frac{R_{M}^{2}}{\alpha_{M-1}}\right)^{-M/4}\\
 & \times & K_{\frac{M}{2}}\left(\sqrt{\left(1-\alpha_{1}\right)\eta_{1}^{2}+\left(\alpha_{1}-\alpha_{2}\right)\eta_{2}^{2}+\cdots+\left(\alpha_{M-2}-\alpha_{M-1}\right)\eta_{M-1}^{2}+\cdots\,+\alpha_{M-1}\eta_{M}^{2}}\right.\quad.\\
 & \times & \left.\sqrt{\frac{R_{1}^{2}}{\left(1-\alpha_{1}\right)}+\frac{R_{2}^{2}}{\left(\alpha_{1}-\alpha_{2}\right)}+\cdots+\frac{R_{M-1}^{2}}{\left(\alpha_{M-2}-\alpha_{M-1}\right)}\,+\frac{R_{M}^{2}}{\alpha_{M-1}}}\right)
\end{array}\label{eq:MtransCompact-Schweber2}
\end{equation}

\section{Unifying the Upper Limits of the Integrals to Give a Third, Numerically-Stable
Integral Representation}

The above forms are numerically stable but, for $M>2$, their utility
is somewhat hampered for analytical reduction since
each succeeding integral has the prior parameter as its upper limit.
One can cast each such integral into one over the interval $[0,1]$
by making a  change of variables to
\begin{equation}
\alpha_{j}\to\alpha_{j-1}\sigma_{j}\label{eq:al_to_sig}
\end{equation}
in sequence from $j=M-1$ down to $j=2$ and by multiplying the set of derivatives of $\alpha_{j}/\alpha_{j-1}$
that defines each new variable $\sigma_{j}$, 
\begin{equation}
\prod_{j=1}^{M-2}\frac{1}{\alpha_{j}}\quad.\label{eq:dervj}
\end{equation}

The first three such are,
where we explicitly put in the shifted coordinates $R_{j}^{2}=x_{1j}^{2}$
(and $\eta_{j}^{2}=\eta_{1j}^{2}$), 
\begin{equation}
\begin{array}{ccc}
\frac{e^{-x_{1}\eta_{1}}}{x_{1}} \frac{e^{-x_{12}\eta_{12}}}{x_{12}} \frac{e^{-x_{13}\eta_{13}}}{x_{13}} & = & \int_{0}^{1}d\alpha_{1}\int_{0}^{1}d\sigma_{2}\frac{\alpha_{1}\left(\left(1-\alpha_{1}\right)\eta_{1}^{2}+\eta_{12}^{2}\alpha_{1}\left(1-\sigma_{2}\right)+\eta_{13}^{2}\alpha_{1}\sigma_{2}\right)^{3/4}}{\sqrt{2}\pi^{3/2}\left(\left(1-\alpha_{1}\right)\alpha_{1}^{2}\left(1-\sigma_{2}\right)\sigma_{2}\right)^{3/2}}\left(\frac{x_{1}^{2}}{1-\alpha_{1}}+\frac{x_{12}^{2}}{\alpha_{1}\left(1-\sigma_{2}\right)}+\frac{x_{13}^{2}}{\alpha_{1}\sigma_{2}}\right)^{-3/4}\\
 & \times & K_{\frac{3}{2}}\left(\sqrt{\frac{x_{1}^{2}}{1-\alpha_{1}}+\frac{x_{12}^{2}}{\alpha_{1}\left(1-\sigma_{2}\right)}+\frac{x_{13}^{2}}{\alpha_{1}\sigma_{2}}}\sqrt{\left(1-\alpha_{1}\right)\eta_{1}^{2}+\eta_{12}^{2}\alpha_{1}\left(1-\sigma_{2}\right)+\eta_{13}^{2}\alpha_{1}\sigma_{2}}\right)\quad,
\end{array}\label{eq:3-sig}
\end{equation}

\begin{equation}
\begin{array}{ccc}
\frac{e^{-x_{1}\eta_{1}}}{x_{1}} \frac{e^{-x_{12}\eta_{12}}}{x_{12}} \frac{e^{-x_{13}\eta_{13}}}{x_{13}} \frac{e^{-x_{14}\eta_{14}}}{x_{14}} & = & \int_{0}^{1}d\alpha_{1}\int_{0}^{1}d\sigma_{2}\int_{0}^{1}d\sigma_{3}\frac{\alpha_{1}^{2}\sigma_{2}\left(\left(1-\alpha_{1}\right)\eta_{1}^{2}+\eta_{12}^{2}\alpha_{1}\left(1-\sigma_{2}\right)+\eta_{13}^{2}\alpha_{1}\sigma_{2}\left(1-\sigma_{3}\right)+\eta_{14}^{2}\alpha_{1}\sigma_{2}\sigma_{3}\right)}{2\pi^{2}\left(\left(1-\alpha_{1}\right)\alpha_{1}^{3}\left(1-\sigma_{2}\right)\sigma_{2}^{2}\left(1-\sigma_{3}\right)\sigma_{3}\right)^{3/2}}\\
 & \times & \hspace{-4.1cm} \left(\frac{x_{1}^{2}}{1-\alpha_{1}}+\frac{x_{12}^{2}}{\alpha_{1}\left(1-\sigma_{2}\right)}+\frac{x_{13}^{2}}{\alpha_{1}\sigma_{2}\left(1-\sigma_{3}\right)}+\frac{x_{14}^{2}}{\alpha_{1}\sigma_{2}\sigma_{3}}\right)^{-1} \\
& \times &   \hspace{-3.0cm} K_{2}\left[\sqrt{\frac{x_{1}^{2}}{1-\alpha_{1}}+\frac{x_{12}^{2}}{\alpha_{1}\left(1-\sigma_{2}\right)}+\frac{x_{13}^{2}}{\alpha_{1}\sigma_{2}\left(1-\sigma_{3}\right)}+\frac{x_{14}^{2}}{\alpha_{1}\sigma_{2}\sigma_{3}}}\right.\\
 & \times &  \left.\sqrt{\left(1-\alpha_{1}\right)\eta_{1}^{2}+\eta_{12}^{2}\alpha_{1}\left(1-\sigma_{2}\right)+\eta_{13}^{2}\alpha_{1}\sigma_{2}\left(1-\sigma_{3}\right)+\eta_{14}^{2}\alpha_{1}\sigma_{2}\sigma_{3}}\right] \quad,
\end{array}\label{eq:4-sig}
\end{equation}
and

\begin{equation}
\begin{array}{ccc}
\frac{e^{-x_{1}\eta_{1}}}{x_{1}} \hspace{-0.35cm} & & \hspace{-5.1cm} \frac{e^{-x_{12}\eta_{12}}}{x_{12}} \frac{e^{-x_{13}\eta_{13}}}{x_{13}} \frac{e^{-x_{14}\eta_{14}}}{x_{14}} \frac{e^{-x_{15}\eta_{15}}}{x_{15}} 
 =  \int_{0}^{1}d\alpha_{1}\int_{0}^{1}d\sigma_{2}\int_{0}^{1}d\sigma_{3}\int_{0}^{1}d\sigma_{4} \\
& \times & \hspace{-2.1cm} \frac{\alpha_{1}^{3}\sigma_{2}^{2}\sigma_{3}\left(\sqrt{\left(1-\alpha_{1}\right)\eta_{1}^{2}+\eta_{12}^{2}\alpha_{1}\left(1-\sigma_{2}\right)+\eta_{13}^{2}\alpha_{1}\sigma_{2}\left(1-\sigma_{3}\right)+\eta_{14}^{2}\alpha_{1}\sigma_{2}\sigma_{3}\left(1-\sigma_{4}\right)+\eta_{15}^{2}\alpha_{1}\sigma_{2}\sigma_{3}\sigma_{4}}\right)^{5/4}}{2\sqrt{2}\pi^{5/2}\left(\left(1-\alpha_{1}\right)\alpha_{1}^{4}\left(1-\sigma_{2}\right)\sigma_{2}^{3}\left(1-\sigma_{3}\right)\sigma_{3}^{2}\left(1-\sigma_{4}\right)\sigma_{4}\right)^{3/2}}\\
 & \times & \hspace{-4.1cm} \left(\frac{x_{1}^{2}}{1-\alpha_{1}}+\frac{x_{12}^{2}}{\alpha_{1}\left(1-\sigma_{2}\right)}+\frac{x_{13}^{2}}{\alpha_{1}\sigma_{2}\left(1-\sigma_{3}\right)}+\frac{x_{14}^{2}}{\alpha_{1}\sigma_{2}\sigma_{3}\left(1-\sigma_{4}\right)}+\frac{x_{15}^{2}}{\alpha_{1}\sigma_{2}\sigma_{3}\sigma_{4}}\right)^{-5/4} \\
& \times & \hspace{-3.8cm} K_{\frac{5}{2}}\left[\sqrt{\frac{x_{1}^{2}}{1-\alpha_{1}}+\frac{x_{12}^{2}}{\alpha_{1}\left(1-\sigma_{2}\right)}+\frac{x_{13}^{2}}{\alpha_{1}\sigma_{2}\left(1-\sigma_{3}\right)}+\frac{x_{14}^{2}}{\alpha_{1}\sigma_{2}\sigma_{3}\left(1-\sigma_{4}\right)}+\frac{x_{15}^{2}}{\alpha_{1}\sigma_{2}\sigma_{3}\sigma_{4}}}\right.\\
& \times & \left.\sqrt{\left(1-\alpha_{1}\right)\eta_{1}^{2}+\eta_{12}^{2}\alpha_{1}\left(1-\sigma_{2}\right)+\eta_{13}^{2}\alpha_{1}\sigma_{2}\left(1-\sigma_{3}\right)+\eta_{14}^{2}\alpha_{1}\sigma_{2}\sigma_{3}\left(1-\sigma_{4}\right)+\eta_{15}^{2}\alpha_{1}\sigma_{2}\sigma_{3}\sigma_{4}}\right] \,.
\end{array}\label{eq:5-sig}
\end{equation}

This can easily be seen to generalize to

\begin{equation}
\begin{array}{ccc}
\frac{e^{-R_{1}\eta_{1}}}{R_{1}}\hspace{-0.4cm} & & \hspace{-6.6cm}\frac{e^{-R_{2}\eta_{2}}}{R_{2}}\cdots\frac{e^{-R_{M}\eta_{M}}}{R_{M}}=\int_{0}^{1}d\alpha_{1}\int_{0}^{1}d\sigma_{2}\int_{0}^{1}d\sigma_{3}\cdots\int_{0}^{1}d\sigma_{M-1} \\
 & \times & \hspace{-4.2cm} {\displaystyle \frac{2^{-\left(M-2\right)/2}\pi^{-M/2}\alpha_{1}^{M-2}\sigma_{2}^{M-3}\sigma_{3}^{M-4}\cdots\sigma_{M-2}}{\left(\left(1-\alpha_{1}\right)\alpha_{1}^{M-1}\left(1-\sigma_{2}\right)\sigma_{2}^{M-2}\left(1-\sigma_{3}\right)\cdots\sigma_{M-2}^{2}\left(1-\sigma_{M-1}\right)\sigma_{M-1}\right)^{3/2}}}\\
 & \times & \hspace{-6.7cm}  \left(\left(1-\alpha_{1}\right)\eta_{1}^{2}+\eta_{2}^{2}\alpha_{1}\left(1-\sigma_{2}\right)+\eta_{3}^{2}\alpha_{1}\sigma_{2}\left(1-\sigma_{3}\right)+\cdots \right.
 \\ 
 & & \left. \hspace{4.3cm} \cdots +\eta_{M-1}^{2}\alpha_{1}\sigma_{2}\sigma_{3}\cdots\sigma_{M-2}\left(1-\sigma_{M-1}\right)  +\eta_{M}^{2}\alpha_{1}\sigma_{2}\sigma_{3}\cdots\sigma_{M-1}\right)^{M/4}\\
 & \times & \hspace{-2.0cm} \left(\frac{R_{1}^{2}}{\left(1-\alpha_{1}\right)}+\frac{R_{2}^{2}}{\alpha_{1}\left(1-\sigma_{2}\right)}+\frac{R_{3}^{2}}{\alpha_{1}\sigma_{2}\left(1-\sigma_{3}\right)}+\cdots+\frac{R_{M-1}^{2}}{\alpha_{1}\sigma_{2}\sigma_{3}\cdots\sigma_{M-2}\left(1-\sigma_{M-1}\right)}+\frac{R_{M}^{2}}{\alpha_{1}\sigma_{2}\sigma_{3}\cdots\sigma_{M-1}}\right)^{-M/4}\\
 & \times & \hspace{-5.9cm}  K_{\frac{M}{2}}\left(
\left[ \left(1-\alpha_{1}\right)\eta_{1}^{2}+\eta_{2}^{2}\alpha_{1}\left(1-\sigma_{2}\right)+\eta_{3}^{2}\alpha_{1}\sigma_{2}\left(1-\sigma_{3}\right)+ \cdots \right.  \right. \\
 &  & \hspace{4.7cm}  \left.  \left. \cdots+\eta_{M-1}^{2} \alpha_{1}\sigma_{2}\sigma_{3}\cdots\sigma_{M-2} \left(1-\sigma_{M-1}\right) + \eta_{M}^{2}\alpha_{1}\sigma_{2}\sigma_{3}\cdots\sigma_{M-1}   \right]^{1/2}  \right.\\
 & \times & \hspace{-2.3cm}  \left.\sqrt{\frac{R_{1}^{2}}{\left(1-\alpha_{1}\right)}+\frac{R_{2}^{2}}{\alpha_{1}\left(1-\sigma_{2}\right)}+\frac{R_{3}^{2}}{\alpha_{1}\sigma_{2}\left(1-\sigma_{3}\right)}+\cdots+\frac{R_{M-1}^{2}}{\alpha_{1}\sigma_{2}\sigma_{3}\cdots\sigma_{M-2}\left(1-\sigma_{M-1}\right)}+\frac{R_{M}^{2}}{\alpha_{1}\sigma_{2}\sigma_{3}\cdots\sigma_{M-1}}}\right) \quad.
\end{array}\label{eq:MtransCompact_sig}
\end{equation}

\noindent Inspection shows that for $M=2$ we indeed obtain eq. (\ref{eq:YYalpha}).

One unusual feature of this integral representation is that, like
that derived in the prior paper, the recursion relationships of Macdonald
functions may be applied to lower (or raise) the indices.

\subsection{Inclusion of plane waves and dipole interactions}

Transition amplitudes sometimes contain plane waves, and these may
be easily included in this integral representation directly in the
$\rho$ version,

\begin{equation}
\begin{array}{ccc}
\frac{e^{-R_{1}\eta_{1}}}{R_{1}}\hspace{-0.3cm} & & \hspace{-4.9cm}\frac{e^{-R_{2}\eta_{2}}}{R_{2}}\cdots\frac{e^{-R_{M}\eta_{M}}}{R_{M}}=\int_{0}^{\infty}d\rho\int_{0}^{1}d\alpha_{1}\int_{0}^{1}d\sigma_{2}\int_{0}^{1}d\sigma_{3}\cdots\int_{0}^{1}d\sigma_{M-1}\\
& \times &  \hspace{-2.9cm} {\displaystyle \frac{2^{-M}\pi^{-M/2}\rho^{-\left(M+2\right)/2}\alpha_{1}^{M-2}\sigma_{2}^{M-3}\sigma_{3}^{M-4}\cdots\sigma_{M-2}}{\left(\left(1-\alpha_{1}\right)\alpha_{1}^{M-1}\left(1-\sigma_{2}\right)\sigma_{2}^{M-2}\left(1-\sigma_{3}\right)\cdots\sigma_{M-2}^{2}\left(1-\sigma_{M-1}\right)\sigma_{M-1}\right)^{3/2}}}\\
 & \times & \hspace{-4.2cm}exp\left(-\rho\left(\left(1-\alpha_{1}\right)\eta_{1}^{2}+\eta_{2}^{2}\alpha_{1}\left(1-\sigma_{2}\right)+\eta_{3}^{2}\alpha_{1}\sigma_{2}\left(1-\sigma_{3}\right)+\cdots \right. \right. \\
 & & \left. \left.  \hspace{3.3cm}  \cdots +\eta_{M-1}^{2}\alpha_{1}\sigma_{2}\sigma_{3}\cdots\sigma_{M-2}\left(1-\sigma_{M-1}\right)+\eta_{M}^{2}\alpha_{1}\sigma_{2}\sigma_{3}\cdots\sigma_{M-1}\right)\right)\\
 & \times &  \hspace{-4.9cm} {\displaystyle exp\left(-
 \left(\frac{R_{1}^{2}}{\left(1-\alpha_{1}\right)}+\frac{R_{2}^{2}}{\alpha_{1}\left(1-\sigma_{2}\right)}+\frac{R_{3}^{2}}{\alpha_{1}\sigma_{2}\left(1-\sigma_{3}\right)}+\cdots \right. \right. }\\
 & & {\displaystyle \left. \left.  \hspace{4.1cm} \cdots +\frac{R_{M-1}^{2}}{\alpha_{1}\sigma_{2}\sigma_{3}\cdots\sigma_{M-2}\left(1-\sigma_{M-1}\right)}+\frac{R_{M}^{2}}{\alpha_{1}\sigma_{2}\sigma_{3}\cdots\sigma_{M-1}}\right)\frac{1}{4\rho}
 \right)
}
\end{array}\label{eq:MtransCompact_sig_rho}
\end{equation}
prior to completing the square, by utilizing an orthogonal transformation (like that for the $k_{j}$)
that reduces the spatial-coordinate quadratic form to diagonal form.
Again, one has invariant determinants for this orthogonal transformation
so that it never needs to actually be explicitly found. As before,
this is followed by a simple translation in $\lbrace\mathbf{x}_{1},\,\mathbf{x}_{2},\cdots,\mathbf{\, x}_{N}\rbrace$
space (with Jacobian = 1).

In the more compact version containing Macdonald functions, one can
simply apply the translation in $\lbrace\mathbf{x}_{1},\,\mathbf{x}_{2},\cdots,\mathbf{\, x}_{N}\rbrace$
space to the plane wave(s) that multiply eqs.  (\ref{eq:MtransCompact-Schweber3}), (\ref{eq:MtransCompact-Schweber2}), and  (\ref{eq:MtransCompact_sig}).

Photoionization transition amplitudes will generally contain dipole
terms $\cos\left(\theta\right)$ that may be transformed into plane
waves via a transformation like $\cos\theta_{12}=-x_{1}^{-1}x_{2}^{-1}\left.\frac{\partial}{\partial Q}e^{-Q{\bf x}_{1}\cdot{\bf x}_{2}}\right|_{Q=0}$,\cite{Stra90a}
giving an integro-differential representation whose inclusion follows
that for other sorts of plane waves.

\section{Utilizing Meijer G-functions to reduce integrals}

The utility of these new integral transformations for large \emph{M}
may well hinge on finding integrals over variables that reside within
square roots as the argument of a Macdonald function. One method for
crafting such untabled integrals is to violate a
two general rules of procedure in analytic reduction of integrals.
The first rule is to use sequential integration whenever possible.
For instance, if one adds a third unshifted Slater orbital to eq.
(\ref{eq:SVxVxy}) and integrates over both variables, one would reasonably
start by transforming only those Slater orbitals that contain $x_{1}$
and integrate over that variable. Next, one integrates the resultant
and the third Slater orbital over $x_{2}$. The result is easily found to be\cite{GR5 p. 358 No. 3.351.3}

\begin{eqnarray}
S_{1}^{\eta_{1}0\eta_{12}0\eta_{2}0}\left(0,0;0,0,0\right) & = & \int d^{3}x_{2}\int d^{3}x_{1}\frac{e^{-\eta_{1}x_{1}}}{x_{1}}\frac{e^{-\eta_{12}x_{12}}}{x_{12}}\frac{e^{-\eta_{2}x_{2}}}{x_{2}}\quad.\label{eq:S_Y1Y12Y2-Seq}\\
 & = & \int_{0}^{\infty}dx_{2}4\pi x_{2}^{2}\frac{4\pi\left(e^{-x_{2}\eta_{12}}-e^{-x_{2}\eta_{1}}\right)}{x_{2}\left(\eta_{1}^{2}-\eta_{12}^{2}\right)}\frac{e^{-\eta_{2}x_{2}}}{x_{2}}=\frac{16\pi^{2}}{\left(\eta_{1}+\eta_{2}\right)\left(\eta_{1}+\eta_{12}\right)\left(\eta_{2}+\eta_{12}\right)}\nonumber 
\end{eqnarray}

There is utility, however, in simultaneously transforming the full
product of Slater orbitals to generate unusual integrals whose values
we know (as above), but whose reduction path may be fraught with difficulty.
If one can find the path for a known integral, this may provide a
path for unknown integrals. And it is clear that the integral
representations of the present paper, like that in the prior paper, are unusual in that they have integration
variables residing within square roots as the arguments of Macdonald
functions. 

The present integral representation appears on the surface
to be less likely to allow for such a reduction because there are
many fewer tabled integrals over the interval $[0,1]$ than there are
over the $[0,\infty]$ interval in the integral representation of the prior work.
We will see in applying this strategy to the above integral that this
concern is not at all the case when one does both coordinate integrals
first. 

We apply the integral representation eq. (\ref{eq:MtransCompact_sig}) to all
three Slater orbitals simultaneously. After completing the square
and changing variables, the integral over $x_{1}^{'2}K_{\frac{3}{2}}\left(\alpha\sqrt{x_{1}^{'2}+z^{2}}\right)/\sqrt{\left(x_{1}^{'2}+z^{2}\right)^{3/2}}$
may be done using \cite{GR5 p. 727 No. 6.596.3 GR7 p. 693} and the
consequent $x_{2}^{2}K_{0}\left(ax_{2}\right)$ integral may also
be done via \cite{GR7 p. 665 No. 6.521.10} , and the second to last
integral is given by \cite{GR5 p. 333 No. 3.194.1 GR7 p. 315} :%

\begin{align}
S_{1}^{\eta_{1}0\eta_{12}0\eta_{2}0}\left(0,0;0,0,0\right) & =\int d^{3}x_{2}\int d^{3}x_{1}\frac{e^{-\eta_{1}x_{1}}}{x_{1}}\frac{e^{-\eta_{12}x_{12}}}{x_{12}}\frac{e^{-\eta_{2}x_{2}}}{x_{2}}\nonumber \\
 & =\int d^{3}x_{2}\int d^{3}x_{1}\int_{0}^{1}d\alpha_{1}\int_{0}^{1}d\sigma_{2}\frac{\alpha_{1}\left(\left(1-\alpha_{1}\right)\eta_{1}^{2}+\eta_{12}^{2}\alpha_{1}\left(1-\sigma_{2}\right)+\eta_{2}^{2}\alpha_{1}\sigma_{2}\right)^{3/4}}{\sqrt{2}\pi^{3/2}\left(\left(1-\alpha_{1}\right)\alpha_{1}^{2}\left(1-\sigma_{2}\right)\sigma_{2}\right)^{3/2}}\nonumber \\
 & \times\frac{K_{\frac{3}{2}}\left(\sqrt{\frac{x_{1}^{2}}{1-\alpha_{1}}+\frac{x_{12}^{2}}{\alpha_{1}\left(1-\sigma_{2}\right)}+\frac{x_{2}^{2}}{\alpha_{1}\sigma_{2}}}\sqrt{\left(1-\alpha_{1}\right)\eta_{1}^{2}+\eta_{12}^{2}\alpha_{1}\left(1-\sigma_{2}\right)+\eta_{2}^{2}\alpha_{1}\sigma_{2}}\right)}{\left(\frac{x_{1}^{2}}{1-\alpha_{1}}+\frac{x_{12}^{2}}{\alpha_{1}\left(1-\sigma_{2}\right)}+\frac{x_{2}^{2}}{\alpha_{1}\sigma_{2}}\right)^{3/4}}\nonumber \\
 & =\int d^{3}x_{2}\int d^{3}x_{1}^{'}\int_{0}^{1}d\alpha_{1}\int_{0}^{1}d\sigma_{2}\frac{\alpha_{1}\left(\left(1-\alpha_{1}\right)\eta_{1}^{2}+\eta_{12}^{2}\alpha_{1}\left(1-\sigma_{2}\right)+\eta_{2}^{2}\alpha_{1}\sigma_{2}\right)^{3/4}}{\sqrt{2}\pi^{3/2}\left(\left(1-\alpha_{1}\right)\alpha_{1}^{2}\left(1-\sigma_{2}\right)\sigma_{2}\right)^{3/2}}\nonumber \\
 & \times\frac{K_{\frac{3}{2}}\left(\sqrt{\frac{x_{1}^{'2}\left(1-\alpha_{1}\sigma_{2}\right)}{\left(1-\alpha_{1}\right)\alpha_{1}\left(1-\sigma_{2}\right)}+\frac{x_{2}^{2}}{\alpha_{1}\sigma_{2}-\alpha_{1}^{2}\sigma_{2}^{2}}}\sqrt{\left(1-\alpha_{1}\right)\eta_{1}^{2}+\eta_{12}^{2}\alpha_{1}\left(1-\sigma_{2}\right)+\eta_{2}^{2}\alpha_{1}\sigma_{2}}\right)}{\left(\frac{x_{1}^{'2}\left(1-\alpha_{1}\sigma_{2}\right)}{\left(1-\alpha_{1}\right)\alpha_{1}\left(1-\sigma_{2}\right)}+\frac{x_{2}^{2}}{\alpha_{1}\sigma_{2}-\alpha_{1}^{2}\sigma_{2}^{2}}\right)^{3/4}}\nonumber \\
 & =\int_{0}^{\infty}dx_{2}\int_{0}^{1}d\alpha_{1}\int_{0}^{1}d\sigma_{2}8\pi x_{2}^{2}K_{0}\left(\frac{x_{2}\sqrt{\left(1-\alpha_{1}\right)\eta_{1}^{2}+\alpha_{1}\eta_{13}^{2}\sigma_{2}+\eta_{12}^{2}\left(\alpha_{1}-\alpha_{1}\sigma_{2}\right)}}{\sqrt{\alpha_{1}}\sqrt{\sigma_{2}}\sqrt{1-\alpha_{1}\sigma_{2}}}\right)\label{eq:YYY-simult-sig}\\
 & \times\frac{\left(\alpha_{1}\eta_{13}^{2}\sigma_{2}+\eta_{12}^{2}\left(\alpha_{1}-\alpha_{1}\sigma_{2}\right)+\left(1-\alpha_{1}\right)\eta_{1}^{2}\right)^{3/4}}{\left(1-\alpha_{1}\right)^{3/4}\alpha_{1}^{5/4}\left(1-\sigma_{2}\right)^{3/4}\sigma_{2}^{3/2}\left(1-\alpha_{1}\sigma_{2}\right)^{3/4}}\nonumber \\
 & \times\left(\frac{\sqrt{1-\alpha_{1}\sigma_{2}}\sqrt{\alpha_{1}\eta_{13}^{2}\sigma_{2}+\eta_{12}^{2}\left(\alpha_{1}-\alpha_{1}\sigma_{2}\right)+\left(1-\alpha_{1}\right)\eta_{1}^{2}}}{\sqrt{1-\alpha_{1}}\sqrt{\alpha_{1}}\sqrt{1-\sigma_{2}}}\right)^{-3/2} \nonumber \\
 & =\int_{0}^{1}d\alpha_{1}\int_{0}^{1}d\sigma_{2}\frac{4\pi^{2}\alpha_{1}}{\left(\alpha_{1}\left(\eta_{12}^{2}\left(1-\sigma_{2}\right)+\eta_{13}^{2}\sigma_{2}-\eta_{1}^{2}\right)+\eta_{1}^{2}\right)^{3/2}}\nonumber \\
 & =\int_{0}^{1}d\sigma_{2}\frac{8\pi^{2}\left(\sqrt{\left(\eta_{13}^{2}-\eta_{12}^{2}\right)\sigma_{2}+\eta_{12}^{2}}-\eta_{1}\right)^{2}}{\left(\left(\eta_{12}^{2}-\eta_{13}^{2}\right)\sigma_{2}+\eta_{1}^{2}-\eta_{12}^{2}\right)^{2}\sqrt{\left(\eta_{13}^{2}-\eta_{12}^{2}\right)\sigma_{2}+\eta_{12}^{2}}}\nonumber \\
 & =\int_{\eta_{12}^{2}}^{\eta_{13}^{2}}\frac{8\pi^{2}}{\sqrt{y}\left(\eta_{13}^{2}-\eta_{12}^{2}\right)\left(\sqrt{y}+\eta_{1}\right)^{2}}\, dy=\int_{\eta_{12}}^{\eta_{13}}\frac{16\pi^{2}}{\left(\eta_{13}^{2}-\eta_{12}^{2}\right)\left(z+\eta_{1}\right)^{2}}\, dz\quad.\nonumber 
\end{align}
Using \cite{GR5 p. 69 No. 2.113.1 GR7 p. 69} in the final integral, above, gives
the final line of eq. (\ref{eq:S_Y1Y12Y2-Seq}), and it does so in
a more straightforward fashion -- using tabled integrals -- than did the
prior integral representation when applied to this problem (that required
the generation of new integrals in the final step using the computer
algebra and calculus program Mathematica).

We will see next that this advantage over the prior method is not
universal. That is, each integral representation, in turn, will shine
brighter on specific problems.

The second general rule of analytical integration is that the easiest
path is to integrate over the coordinate variables first. The doubly
contrary approach of the prior paper generated a set of integrals
that might be of utility for future researchers. We showed therein
that integrating (over the interval $[0,\infty]$) over variables
that reside within square roots as the argument of a Macdonald function
can be done if we rewrite the Macdonald function in terms of a Meijer
G-function. In the case of a product of three Slater orbitals, after
completing the square in the variable common to all three and then
integrating over the result one has the following result to be integrated
over: \cite{PBM3 p. 665 No. 8.4.23.1} 
\begin{equation}
\frac{1}{\zeta_{2}^{3/2}}K_{0}\left(2\frac{x_{2}\sqrt{\zeta_{1}+\zeta_{2}+1}\eta_{2}\sqrt{\frac{\zeta_{1}\eta_{12}^{2}+\eta_{1}^{2}}{4\eta_{2}^{2}}+\frac{\zeta_{2}}{4}}}{\sqrt{\zeta_{1}+1}\sqrt{\zeta_{2}}}\right)=\frac{1}{2}\frac{1}{\zeta_{2}^{3/2}}G_{0,2}^{2,0}\left(\frac{x_{2}^{2}\left(\zeta_{1}+\zeta_{2}+1\right)\eta_{2}^{2}\left(\zeta_{2}+\frac{\eta_{1}^{2}+\zeta_{1}\eta_{12}^{2}}{\eta_{2}^{2}}\right)}{\left(\zeta_{1}+1\right)\zeta_{2}}|\begin{array}{c}
0,0\end{array}\right)\quad, \label{eq:K0 as G2002}
\end{equation}
 for which there is but one tabled integral\cite{PBM3 p. 349 No. 2.24.2.9} that has roughly the right
form (with $\zeta_{2}=x$),

\begin{align}
\int_{0}^{\infty}dx\, x^{\alpha-1}\left(ax^{2}+bx+c\right)^{\frac{3}{2}-\alpha} & G_{0,2}^{2,0}\left(\frac{ax^{2}+bx+c}{x}|\begin{array}{c}
\nu,-\nu\end{array}\right)\nonumber \\
 & =\frac{\sqrt{\pi}G_{1,3}^{3,0}\left(b+2\sqrt{a}\sqrt{c}|\begin{array}{c}
\frac{3}{2}\\
0,-\alpha-\nu+3,-\alpha+\nu+3
\end{array}\right)}{2a^{3/2}}\quad.\label{eq:PBM3 p. 349 No. 2.24.2.9}\\
 & +\frac{\sqrt{\pi}\sqrt{c}G_{1,3}^{3,0}\left(b+2\sqrt{a}\sqrt{c}|\begin{array}{c}
\frac{1}{2}\\
0,-\alpha-\nu+2,-\alpha+\nu+2
\end{array}\right)}{a}\nonumber 
\end{align}
A modification was required since inserting $\alpha=\frac{3}{2}$
to remove the polynomial multiplying the G-function in the integrand
leaves us with the wrong power of \emph{x}. One may, however, take
derivatives with respect to \emph{c} of the integrand and resultant,
with $\nu=1/2$ in combination with $\nu=0$, to show that%
\footnote{In eq. (\ref{eq:K0/x^3/2}) and following, we explicitly include an
alternative expression for K as the hypergeometric U function,\cite{functions.wolfram.com/03.04.26.0003.01}
while expressions in related functions may be found at \cite{functions.wolfram.com/03.04.26.0002.01},
\cite{functions.wolfram.com/03.04.26.0005.01}, and \cite{GR5 p. 1090 No.9.235.2}.
None of these seem to have tabled integrals with arguments as complicated
as \cite{PBM3 p. 349 No. 2.24.2.9}.%
}
%

\begin{align}
\int_{0}^{\infty}dx\,\frac{1}{x^{3/2}}K_{0}\left(2\sqrt{\frac{ax^{2}+bx+c}{x}}\right) & =\int_{0}^{\infty}dx\,\frac{1}{x^{3/2}}\sqrt{\pi}e^{-2\sqrt{\frac{ax^{2}+bx+c}{x}}}U\left(\frac{1}{2},1,4\sqrt{\frac{ax^{2}+bx+c}{x}}\right)\nonumber \\
 & =\int_{0}^{\infty}dx\,\frac{1}{2x^{3/2}}G_{0,2}^{2,0}\left(\frac{\text{ax}^{2}+\text{b x}+c}{x}|\begin{array}{c}
0,0\end{array}\right)\nonumber \\
 & =\frac{\pi e^{-2\sqrt{2\sqrt{a}\sqrt{c}+b}}}{2\sqrt{c}\sqrt{2\sqrt{a}\sqrt{c}+b}}-\frac{\sqrt{\pi}G_{1,3}^{2,1}\left(b+2\sqrt{a}\sqrt{c}|\begin{array}{c}
-\frac{3}{2}\\
-\frac{1}{2},0,-\frac{1}{2}
\end{array}\right)}{2\sqrt{c}}\quad,\label{eq:K0/x^3/2}\\
 & =\frac{\pi e^{-2\sqrt{2\sqrt{a}\sqrt{c}+b}}}{2\sqrt{c}}\nonumber \\
 & \Rightarrow\frac{\pi e^{-2\left(\frac{2\sqrt{a}\sqrt{c}}{x_{2}\eta_{2}}+\frac{x_{2}\eta_{2}}{2}\right)}}{2\sqrt{c}}\nonumber 
\end{align}
where the reduction of the Meijer G-function in the third line is
from \cite{functions.wolfram.com/07.34.03.0727.01} and the last step
\begin{equation}
\sqrt{2\sqrt{a}\sqrt{c}+b}\Rightarrow\frac{2\sqrt{a}\sqrt{c}}{x_{2}\eta_{2}}+\frac{x_{2}\eta_{2}}{2}\label{eq:undoing sqrt}
\end{equation}
 holds for a number of cases akin to the present one in which 
\begin{equation}
\left\{ a\to\frac{x_{2}^{2}\eta_{2}^{2}}{4\left(\zeta_{1}+1\right)},\, b\to\frac{x_{2}^{2}\eta_{2}^{2}}{4\left(\zeta_{1}+1\right)}\left(\frac{\zeta_{1}\eta_{12}^{2}+\eta_{1}^{2}}{\eta_{2}^{2}}+\zeta_{1}+1\right),\, c\to\frac{1}{4}x_{2}^{2}\left(\zeta_{1}\eta_{12}^{2}+\eta_{1}^{2}\right)\right\} \quad.\label{eq:abc->}
\end{equation}

The integral representations of the current paper use integrals over
the interval $[0,1]$ rather than over $[0,\infty]$, and the only
tabled integral over a similar G-function on $[0,1]$ we found\cite{PBM3 p. 349 No. 2.24.2.7} 

\begin{equation}
\int_{0}^{1}dx\,\frac{x^{\alpha-1}\left(1-x\right)^{\beta-1}}{\left(1+ax+b\left(1-x\right)\right)^{\alpha+\beta}}G_{p,q}^{m,n}\left(\frac{x^{\ell}\left(1-x\right)^{k}}{\left(1+ax+b\left(1-x\right)\right)^{\ell+k}}\left|\begin{array}{c}
\left(a_{p}\right)\\
\left(b_{q}\right)
\end{array}\right.\right)\label{eq:G_over[0,1]}
\end{equation}
 has a markedly different integrand than the square of the argument
of the Macdonald function 

\begin{equation}
\frac{8\pi x_{2}^{2}}{\sqrt{\alpha_{1}}\sigma_{2}^{3/2}\left(1-\alpha_{1}\sigma_{2}\right)^{3/2}}K_{0}\left(\frac{x_{2}\sqrt{\left(1-\alpha_{1}\right)\eta_{1}^{2}+\alpha_{1}\eta_{13}^{2}\sigma_{2}+\eta_{12}^{2}\left(\alpha_{1}-\alpha_{1}\sigma_{2}\right)}}{\sqrt{\alpha_{1}}\sqrt{\sigma_{2}}\sqrt{1-\alpha_{1}\sigma_{2}}}\right)\label{eq:K_0_alph_sig}
\end{equation}
to which the fourth equality in (\ref{eq:YYY-simult-sig}) reduces,
so it is useless for the present problem. 

Schweber's third parametrization gives no better result. All of this
provided motivation for the integral representation we derive
in the next section.

\section{An Ungainly but Useful Bridge}

There is an obscure integral\cite{GH p. 176 No. 421.8} that will serve as an integral representation for a product of denominators,

\begin{equation}
\begin{array}{ccc}
\frac{B(\kappa,\lambda-\kappa)}{c^{\kappa}d^{\lambda-\kappa}} & = & \int_{0}^{1}dx\,\left(\frac{x^{\kappa-1}}{(cx+d)^{\lambda}}+\frac{x^{-\kappa+\lambda-1}}{(c+dx)^{\lambda}}\right)=\int_{1}^{\infty}dx\,\left(\frac{x^{\kappa-1}}{(cx+d)^{\lambda}}+\frac{x^{-\kappa+\lambda-1}}{(c+dx)^{\lambda}}\right)\\
 & = & \frac{1}{2}\int_{0}^{\infty}dx\,\left(\frac{x^{\kappa-1}}{(cx+d)^{\lambda}}+\frac{x^{-\kappa+\lambda-1}}{(c+dx)^{\lambda}}\right)\quad,
\end{array}\label{eq:GH p. 176 No. 421.8}
\end{equation}
that has the very useful property of relating integrals over $[0,\infty]$
to integrals over $[0,1]$ and $[1,\infty]$. One might hope to extend
integrals like \cite{PBM3 p. 349 No. 2.24.2.9} to integrals over
$[0,1]$ or $[1,\infty]$ if the integral over $[0,\infty]$  could be
found.

Its ungainliness is revealed in the process of extending it from a
pair of denominators
\begin{equation}
\frac{1}{a_{1}a_{2}}=\int_{0}^{1}d\sigma_{1}\left(\frac{1}{\left(a_{1}+a_{2}\sigma_{1}\right)^{2}}+\frac{1}{\left(a_{1}\sigma_{1}+a_{2}\right)^{2}}\right)\label{eq:pair}
\end{equation}
to triplets by iteration 

\begin{eqnarray}
\frac{1}{a_{1}a_{2}a_{3}} & = &\int_{0}^{1}d\sigma_{1}\int_{0}^{1}d\sigma_{2}2\left(\frac{1}{\left(a_{1}\sigma_{1}+a_{3}\sigma_{2}+a_{2}\right)^{3}}+\frac{1}{\left(a_{2}\sigma_{1}+a_{3}\sigma_{2}+a_{1}\right)^{3}} \right. \nonumber \\
& + & \left.\frac{\sigma_{2}} {\left(\sigma_{2}\left(a_{1}\sigma_{1}+a_{2}\right)+a_{3}\right)^{3}}+\frac{\sigma_{2}}{\left(\sigma_{2}\left(a_{2}\sigma_{1}+a_{1}\right)+a_{3}\right)^{3}}\right)\label{eq:pair-1}
\end{eqnarray}
and beyond: at each step, the number of terms doubles so a general-M
version is difficult to imagine. 

To derive an integral representation for a product of two or three
Slater orbitals, one simply follows the procedure laid out in Section
4, above, for each of the two or four terms, but the determinants
are different in this new case and are different from term to term.
The final forms are%

\begin{align}
\begin{array}{ccc}
{\displaystyle \frac{e^{-\eta_{1}x_{1}}}{x_{1}}\frac{e^{-\eta_{12}x_{12}}}{x_{12}}} & = & \int_{0}^{1}d\sigma_{1}\int_{0}^{\infty}d\rho\frac{\exp\left(-\frac{\frac{x_{1}^{2}}{\sigma_{1}}+x_{12}^{2}}{4\rho}-\rho\left(\eta_{1}^{2}\sigma_{1}+\eta_{12}^{2}\right)\right)+\exp\left(-\frac{\frac{x_{12}^{2}}{\sigma_{1}}+x_{1}^{2}}{4\rho}-\rho\left(\eta_{12}^{2}\sigma_{1}+\eta_{1}^{2}\right)\right)}{4\pi\rho^{2}\sigma_{1}^{3/2}}\\
 & = & \int_{0}^{1}d\sigma_{1}\left(\frac{\sqrt{\eta_{1}^{2}\sigma_{1}+\eta_{12}^{2}}K_{1}\left(\sqrt{\frac{x_{1}^{2}}{\sigma_{1}}+x_{12}^{2}}\sqrt{\sigma_{1}\eta_{1}^{2}+\eta_{12}^{2}}\right)}{\pi\sigma_{1}^{3/2}\sqrt{\frac{x_{1}^{2}}{\sigma_{1}}+x_{12}^{2}}}+\frac{\sqrt{\eta_{12}^{2}\sigma_{1}+\eta_{1}^{2}}K_{1}\left(\sqrt{x_{1}^{2}+\frac{x_{12}^{2}}{\sigma_{1}}}\sqrt{\eta_{1}^{2}+\eta_{12}^{2}\sigma_{1}}\right)}{\pi\sigma_{1}^{3/2}\sqrt{\frac{x_{12}^{2}}{\sigma_{1}}+x_{1}^{2}}}\right)
\end{array}\label{eq:YYsig}
\end{align}

and

\begin{align}
\begin{array}{ccc}
{\displaystyle \frac{e^{-\eta_{1}x_{1}}}{x_{1}}\frac{e^{-\eta_{12}x_{12}}}{x_{12}}\frac{e^{-\eta_{13}x_{13}}}{x_{13}}} & = & \hspace{-4.2cm} \int_{0}^{1}d\sigma_{1}\int_{0}^{1}d\sigma_{2}\int_{0}^{\infty}d\rho{\displaystyle \frac{1}{8\pi^{3/2}\rho^{5/2}\left(\sigma_{1}\sigma_{2}\right){}^{3/2}}}\\
 & + & \left(\exp\left(-\frac{\frac{x_{1}^{2}}{\sigma_{1}}+\frac{x_{13}^{2}}{\sigma_{2}}+x_{12}^{2}}{4\rho}-\rho\left(\eta_{1}^{2}\sigma_{1}+\eta_{13}^{2}\sigma_{2}+\eta_{12}^{2}\right)\right)\right.\\
 & + & \exp\left(-\frac{\frac{x_{12}^{2}}{\sigma_{1}}+\frac{x_{13}^{2}}{\sigma_{2}}+x_{1}^{2}}{4\rho}-\rho\left(\eta_{12}^{2}\sigma_{1}+\eta_{13}^{2}\sigma_{2}+\eta_{1}^{2}\right)\right)\\
 & + & \sigma_{2}\exp\left(-\frac{\frac{x_{1}^{2}}{\sigma_{1}\sigma_{2}}+\frac{x_{12}^{2}}{\sigma_{2}}+x_{13}^{2}}{4\rho}-\rho\left(\eta_{1}^{2}\sigma_{1}\sigma_{2}+\eta_{12}^{2}\sigma_{2}+\eta_{13}^{2}\right)\right)\\
 & + & \left.\sigma_{2}\exp\left(-\frac{\frac{x_{1}^{2}}{\sigma_{2}}+\frac{x_{12}^{2}}{\sigma_{1}\sigma_{2}}+x_{13}^{2}}{4\rho}-\rho\left(\eta_{1}^{2}\sigma_{2}+\eta_{12}^{2}\sigma_{1}\sigma_{2}+\eta_{13}^{2}\right)\right)\right)\\
 & = & \int_{0}^{1}d\sigma_{1}\int_{0}^{1}d\sigma_{2}\left(\frac{\left(\eta_{1}^{2}\sigma_{1}+\eta_{13}^{2}\sigma_{2}+\eta_{12}^{2}\right){}^{3/4}K_{\frac{3}{2}}\left(\sqrt{\frac{x_{1}^{2}}{\sigma_{1}}+x_{12}^{2}+\frac{x_{13}^{2}}{\sigma_{2}}}\sqrt{\sigma_{1}\eta_{1}^{2}+\eta_{12}^{2}+\eta_{13}^{2}\sigma_{2}}\right)}{\sqrt{2}\pi^{3/2}\left(\sigma_{1}\sigma_{2}\right){}^{3/2}\left(\frac{x_{1}^{2}}{\sigma_{1}}+\frac{x_{13}^{2}}{\sigma_{2}}+x_{12}^{2}\right){}^{3/4}}\right.\\
 & + & \frac{\left(\eta_{12}^{2}\sigma_{1}+\eta_{13}^{2}\sigma_{2}+\eta_{1}^{2}\right){}^{3/4}K_{\frac{3}{2}}\left(\sqrt{x_{1}^{2}+\frac{x_{12}^{2}}{\sigma_{1}}+\frac{x_{13}^{2}}{\sigma_{2}}}\sqrt{\eta_{1}^{2}+\eta_{12}^{2}\sigma_{1}+\eta_{13}^{2}\sigma_{2}}\right)}{\sqrt{2}\pi^{3/2}\left(\sigma_{1}\sigma_{2}\right){}^{3/2}\left(\frac{x_{12}^{2}}{\sigma_{1}}+\frac{x_{13}^{2}}{\sigma_{2}}+x_{1}^{2}\right){}^{3/4}}\\
 &  & \frac{\left(\eta_{1}^{2}\sigma_{1}\sigma_{2}+\eta_{12}^{2}\sigma_{2}+\eta_{13}^{2}\right){}^{3/4}K_{\frac{3}{2}}\left(\sqrt{\frac{x_{1}^{2}}{\sigma_{1}\sigma_{2}}+x_{13}^{2}+\frac{x_{12}^{2}}{\sigma_{2}}}\sqrt{\sigma_{1}\sigma_{2}\eta_{1}^{2}+\eta_{13}^{2}+\eta_{12}^{2}\sigma_{2}}\right)}{\sqrt{2}\pi^{3/2}\sigma_{1}^{3/2}\sigma_{2}^{2}\left(\frac{x_{1}^{2}}{\sigma_{1}\sigma_{2}}+\frac{x_{12}^{2}}{\sigma_{2}}+x_{13}^{2}\right){}^{3/4}}\\
 & + & \left.\frac{\left(\eta_{1}^{2}\sigma_{2}+\eta_{12}^{2}\sigma_{1}\sigma_{2}+\eta_{13}^{2}\right){}^{3/4}K_{\frac{3}{2}}\left(\sqrt{\frac{x_{1}^{2}}{\sigma_{2}}+x_{13}^{2}+\frac{x_{12}^{2}}{\sigma_{1}\sigma_{2}}}\sqrt{\sigma_{2}\eta_{1}^{2}+\eta_{13}^{2}+\eta_{12}^{2}\sigma_{1}\sigma_{2}}\right)}{\sqrt{2}\pi^{3/2}\sigma_{1}^{3/2}\sigma_{2}^{2}\left(\frac{x_{1}^{2}}{\sigma_{2}}+\frac{x_{12}^{2}}{\sigma_{1}\sigma_{2}}+x_{13}^{2}\right){}^{3/4}}\right)
\end{array} & \quad,\label{eq:YYYsig}
\end{align}
with the equality also holding for integrals over $[1,\infty]$. It
also holds over $[0,\infty]$ if one multiplies the right-hand side by
 $\frac{1}{2}\frac{1}{2}$. In the case of $M\geq3$, one may even mix these three intervals
among the integrals present.

We (simultaneously) apply the above integral representation (\ref{eq:YYYsig})
to all three Slater orbitals in the first line of (\ref{eq:S_Y1Y12Y2-Seq}),
whose second line gives $117.4952904891590$ when we abitrarily set
parameters to $\left\{ \eta_{1}\to0.3,\eta_{12}\to0.5,\right. $ $\left.\eta_{13}\to0.9\right\} $.
After completing the square and changing variables, using \cite{GR5 p. 727 No. 6.596.3 GR7 p. 693} one may do the integral over

\noindent
$x_{1}^{'2}K_{\frac{3}{2}}\left(\alpha\sqrt{x_{1}^{'2}+z^{2}}\right)/\sqrt{\left(x_{1}^{'2}+z^{2}\right)^{3/2}}$
. %

\begin{align}
& S_{1}^{\eta_{1}0\eta_{12}0\eta_{2}0}\left(0,0;0,0,0\right)  =\int d^{3}x_{2}\int d^{3}x_{1}\frac{e^{-\eta_{1}x_{1}}}{x_{1}}\frac{e^{-\eta_{12}x_{12}}}{x_{12}}\frac{e^{-\eta_{2}x_{2}}}{x_{2}}=\int d^{3}x_{2}\int d^{3}x'_{1}\int_{0}^{1}d\sigma_{1}\int_{0}^{1}d\sigma_{2}\nonumber \\
 & \times\left(\frac{\left(\eta_{1}^{2}\sigma_{1}+\eta_{13}^{2}\sigma_{2}+\eta_{12}^{2}\right){}^{3/4}K_{\frac{3}{2}}\left(\sqrt{\frac{x_{1}^{'2}\left(\sigma_{1}+1\right)}{\sigma_{1}}+\frac{x_{2}^{2}\left(\sigma_{1}+\sigma_{2}+1\right)}{\left(\sigma_{1}+1\right)\sigma_{2}}}\sqrt{\sigma_{1}\eta_{1}^{2}+\eta_{12}^{2}+\eta_{13}^{2}\sigma_{2}}\right)}{\sqrt{2}\pi^{3/2}\left(\sigma_{1}\sigma_{2}\right){}^{3/2}\left(\frac{x_{1}^{'2}\left(\sigma_{1}+1\right)}{\sigma_{1}}+\frac{x_{2}^{2}\left(\sigma_{1}+\sigma_{2}+1\right)}{\left(\sigma_{1}+1\right)\sigma_{2}}\right){}^{3/4}}\right.\nonumber \\
 & +\frac{\left(\eta_{12}^{2}\sigma_{1}+\eta_{13}^{2}\sigma_{2}+\eta_{1}^{2}\right){}^{3/4}K_{\frac{3}{2}}\left(\sqrt{\frac{x_{1}^{'2}\left(\sigma_{1}+1\right)}{\sigma_{1}}+\frac{x_{2}^{2}\left(\sigma_{1}+\sigma_{2}+1\right)}{\left(\sigma_{1}+1\right)\sigma_{2}}}\sqrt{\eta_{1}^{2}+\eta_{12}^{2}\sigma_{1}+\eta_{13}^{2}\sigma_{2}}\right)}{\sqrt{2}\pi^{3/2}\left(\sigma_{1}\sigma_{2}\right){}^{3/2}\left(\frac{x_{1}^{'2}\left(\sigma_{1}+1\right)}{\sigma_{1}}+\frac{x_{2}^{2}\left(\sigma_{1}+\sigma_{2}+1\right)}{\left(\sigma_{1}+1\right)\sigma_{2}}\right){}^{3/4}}\nonumber \\
 & +\frac{\left(\eta_{1}^{2}\sigma_{1}\sigma_{2}+\eta_{12}^{2}\sigma_{2}+\eta_{13}^{2}\right){}^{3/4}K_{\frac{3}{2}}\left(\sqrt{\frac{x_{1}^{'2}\left(\sigma_{1}+1\right)}{\sigma_{1}\sigma_{2}}+\frac{x_{2}^{2}\left(\left(\sigma_{1}+1\right)\sigma_{2}+1\right)}{\left(\sigma_{1}+1\right)\sigma_{2}}}\sqrt{\sigma_{1}\sigma_{2}\eta_{1}^{2}+\eta_{13}^{2}+\eta_{12}^{2}\sigma_{2}}\right)}{\sqrt{2}\pi^{3/2}\sigma_{1}^{3/2}\sigma_{2}^{2}\left(\frac{x_{1}^{'2}\left(\sigma_{1}+1\right)}{\sigma_{1}\sigma_{2}}+\frac{x_{2}^{2}\left(\left(\sigma_{1}+1\right)\sigma_{2}+1\right)}{\left(\sigma_{1}+1\right)\sigma_{2}}\right){}^{3/4}}\nonumber \\
 & +\left.\frac{\left(\eta_{1}^{2}\sigma_{2}+\eta_{12}^{2}\sigma_{1}\sigma_{2}+\eta_{13}^{2}\right){}^{3/4}K_{\frac{3}{2}}\left(\sqrt{\frac{x_{1}^{'2}\left(\sigma_{1}+1\right)}{\sigma_{1}\sigma_{2}}+\frac{x_{2}^{2}\left(\left(\sigma_{1}+1\right)\sigma_{2}+1\right)}{\left(\sigma_{1}+1\right)\sigma_{2}}}\sqrt{\sigma_{2}\eta_{1}^{2}+\eta_{13}^{2}+\eta_{12}^{2}\sigma_{1}\sigma_{2}}\right)}{\sqrt{2}\pi^{3/2}\sigma_{1}^{3/2}\sigma_{2}^{2}\left(\frac{x_{1}^{'2}\left(\sigma_{1}+1\right)}{\sigma_{1}\sigma_{2}}+\frac{x_{2}^{2}\left(\left(\sigma_{1}+1\right)\sigma_{2}+1\right)}{\left(\sigma_{1}+1\right)\sigma_{2}}\right){}^{3/4}}\right)\nonumber \\
 & =\int d^{3}x_{2}\int_{0}^{1}d\sigma_{1}\int_{0}^{1}d\sigma_{2}\left(\frac{2K_{0}\left(\frac{x_{2}\sqrt{\sigma_{1}+\sigma_{2}+1}\sqrt{\sigma_{1}\eta_{1}^{2}+\eta_{12}^{2}+\eta_{13}^{2}\sigma_{2}}}{\sqrt{\sigma_{1}+1}\sqrt{\sigma_{2}}}\right)}{\left(\sigma_{1}+1\right){}^{3/2}\sigma_{2}^{3/2}}\right.\nonumber \\
 & +\frac{2K_{0}\left(\frac{x_{2}\sqrt{\sigma_{1}+\sigma_{2}+1}\sqrt{\eta_{1}^{2}+\eta_{12}^{2}\sigma_{1}+\eta_{13}^{2}\sigma_{2}}}{\sqrt{\sigma_{1}+1}\sqrt{\sigma_{2}}}\right)}{\left(\sigma_{1}+1\right){}^{3/2}\sigma_{2}^{3/2}} \nonumber \\
 & +\left.\frac{2K_{0}\left(\frac{x_{2}\sqrt{\sigma_{1}\sigma_{2}\eta_{1}^{2}+\eta_{13}^{2}+\eta_{12}^{2}\sigma_{2}}\sqrt{\left(\sigma_{1}+1\right)\sigma_{2}+1}}{\sqrt{\sigma_{1}+1}\sqrt{\sigma_{2}}}\right)}{\left(\sigma_{1}+1\right){}^{3/2}\sqrt{\sigma_{2}}}+\frac{2K_{0}\left(\frac{x_{2}\sqrt{\sigma_{2}\eta_{1}^{2}+\eta_{13}^{2}+\eta_{12}^{2}\sigma_{1}\sigma_{2}}\sqrt{\left(\sigma_{1}+1\right)\sigma_{2}+1}}{\sqrt{\sigma_{1}+1}\sqrt{\sigma_{2}}}\right)}{\left(\sigma_{1}+1\right){}^{3/2}\sqrt{\sigma_{2}}}\right)\,. \label{eq:Y1Y12Y2Simultaneous_sig}
\end{align}

For the values of $\left\{ \eta_{1}\to0.3,\eta_{12}\to0.5,\eta_{13}\to0.9\right\} $ used
above, these four terms numerically integrate to 

\begin{equation}
S_{1}^{0.3,0,0.5,0,0.9,0}\left(0,0;0,0,0\right)=\left(39.2072+61.8386+7.89946+8.55004\right)=117.49528665800858\quad.\label{eq:S3_1}
\end{equation}
On the other hand, if we change the integration limits to $\left[1,\infty\right]$
they yield

\begin{equation}
S_{1}^{0.3,0,0.5,0,0.9,0}\left(0,0;0,0,0\right)=\left(31.4147+22.115+38.9735+24.9916\right)=117.49480820934275\quad.\label{eq:S3_1-1}
\end{equation}
If we now change the integration limits to $\left[0,\infty\right]$
and multiply by $\frac{1}{2}\frac{1}{2}$, these terms are

\begin{eqnarray}
S_{1}^{0.3,0,0.5,0,0.9,0}\left(0,0;0,0,0\right) & = & \left(29.373735823279375+29.37376872163303+29.373735823279382 \right. \\
&  &  \left. \hspace{3.3cm} + \,  29.37376872163304\right)=117.49500908982483\quad. \nonumber \label{eq:S3inf}
\end{eqnarray}

It should not be surprising that changing the limits of integration
will change the value of an integral, nor that different integrands
will yield different results under changed limits of integration.
What is remarkable is that the four different integrands in each case
compensate for each other under such a change of integration limits
so as to produce the same sum (to six-digit accuracy in this numerical
integration) for all three limit sets.%
\footnote{Such compensation is inherent in the integral representation (\ref{eq:pair})
we used as the stepping stone to the Slater orbital version expressed
in equations (\ref{eq:YYsig}) and (\ref{eq:YYYsig}).%
} This compensation is also a harbinger that the $\left[0,\infty\right]$
form of the present set of integral representations will indeed act
as a bridge to analytical results for \emph{sums} of integrals of
Meijer G-function over the intervals $\left[0,1\right]$ and$\left[1,\infty\right]$.
It is also notable that the four terms contribute almost equally when
the limits are set to $\left[0,\infty\right]$, with the first and
third terms the same to 15 digits -- as are the second and fourth --
though the first and second terms differ in the seveth digit.

\section{Analytical Results for Sums of Integrals of \textmd{\normalsize{} $f\left(x\right)G_{0,2}^{2,0}\left(\frac{\text{ax}^{2}+\text{b x}+c}{x}|\protect\begin{array}{c}
0,0\protect\end{array}\right)$} over the intervals \textmd{\normalsize $\left[0,1\right]$} and\textmd{\normalsize $\left[1,\infty\right]$}.}

Rather than doing the $x_{2}^{2}K_{0}\left(ax_{2}\right)$ integrals
via \cite{GR7 p. 665 No. 6.521.10}, we will attempt to do the $\sigma_{2}$
integrals first as a means to generate some useful integrals for future
researchers. We begin by recasting the Macdonald functions of eq. (\ref{eq:Y1Y12Y2Simultaneous_sig})
as Meijer G-functions using (\ref{eq:K0 as G2002}), \cite{PBM3 p. 665 No. 8.4.23.1}
with each of the four terms having different values for

\begin{equation}
\begin{array}{c}
\left\{ a\to\frac{x_{2}^{2}\eta_{13}^{2}}{4\left(\sigma_{1}+1\right)},\, b\to\frac{x_{2}^{2}\left(\eta_{1}^{2}\sigma_{1}+\eta_{13}^{2}\left(\sigma_{1}+1\right)+\eta_{12}^{2}\right)}{4\left(\sigma_{1}+1\right)},\, c\to\frac{1}{4}x_{2}^{2}\left(\eta_{1}^{2}\sigma_{1}+\eta_{12}^{2}\right)\right\} \\
\left\{ a_{2}\to\frac{x_{2}^{2}\eta_{13}^{2}}{4\left(\sigma_{1}+1\right)},\, b_{2}\to\frac{x_{2}^{2}\left(\eta_{12}^{2}\sigma_{1}+\eta_{13}^{2}\left(\sigma_{1}+1\right)+\eta_{1}^{2}\right)}{4\left(\sigma_{1}+1\right)},\, c_{2}\to\frac{1}{4}x_{2}^{2}\left(\eta_{12}^{2}\sigma_{1}+\eta_{1}^{2}\right)\right\} \\
\left\{ a_{3}\to\frac{1}{4}x_{2}^{2}\left(\eta_{1}^{2}\sigma_{1}+\eta_{12}^{2}\right),\, b_{3}\to\frac{x_{2}^{2}\left(\eta_{1}^{2}\sigma_{1}+\eta_{13}^{2}\left(\sigma_{1}+1\right)+\eta_{12}^{2}\right)}{4\left(\sigma_{1}+1\right)},\, c_{3}\to\frac{x_{2}^{2}\eta_{13}^{2}}{4\left(\sigma_{1}+1\right)}\right\} \\
\left\{ a_{4}\to\frac{1}{4}x_{2}^{2}\left(\eta_{12}^{2}\sigma_{1}+\eta_{1}^{2}\right),\, b_{4}\to\frac{x_{2}^{2}\left(\eta_{12}^{2}\sigma_{1}+\eta_{13}^{2}\left(\sigma_{1}+1\right)+\eta_{1}^{2}\right)}{4\left(\sigma_{1}+1\right)},\, c_{4}\to\frac{x_{2}^{2}\eta_{13}^{2}}{4\left(\sigma_{1}+1\right)}\right\} 
\end{array}\quad.\label{eq:abc->1-4}
\end{equation}

We switch to the form with integrals over $\left[0,\infty\right]$
so that we can consider just the first term of (\ref{eq:Y1Y12Y2Simultaneous_sig})
that is precisely the integral found in the prior paper eq. (\ref{eq:K0/x^3/2})
except that in the present case, with $\left\{ a\to\frac{x_{2}^{2}\eta_{13}^{2}}{4\left(\sigma_{1}+1\right)},\, b\to\frac{x_{2}^{2}\left(\eta_{1}^{2}\sigma_{1}+\eta_{13}^{2}\left(\sigma_{1}+1\right)+\eta_{12}^{2}\right)}{4\left(\sigma_{1}+1\right)},\, c\to\frac{1}{4}x_{2}^{2}\left(\eta_{1}^{2}\sigma_{1}+\eta_{12}^{2}\right)\right\} $,
the last term simplifies differently:

\begin{equation}
\sqrt{2\sqrt{a}\sqrt{c}+b}\Rightarrow\sqrt{a\left(\sigma_{1}+1\right)}+2\sqrt{\frac{c}{\sigma_{1}+1}}\quad.\label{eq:undoing sqrt-1}
\end{equation}
Then

\begin{equation}
\begin{array}{ccc}
\int d^{3}x_{2}\frac{1}{2}\frac{1}{2}\int_{0}^{\infty}d\sigma_{1}\int_{0}^{\infty}d\sigma_{2}\frac{2}{\left(\sigma_{1}+1\right)^{3/2}\sigma_{2}^{3/2}}K_{0}\left(\frac{2\sqrt{\frac{\eta_{13}^{2}\sigma_{2}^{2}x_{2}^{2}}{4\left(\sigma_{1}+1\right)}+\frac{1}{4}\left(\sigma_{1}\eta_{1}^{2}+\eta_{12}^{2}\right)x_{2}^{2}+\frac{\left(\sigma_{1}\eta_{1}^{2}+\eta_{12}^{2}+\eta_{13}^{2}\left(\sigma_{1}+1\right)\right)\sigma_{2}x_{2}^{2}}{4\left(\sigma_{1}+1\right)}}}{\sqrt{\sigma_{2}}}\right) & =\\
\int d^{3}x_{2}\frac{1}{2}\frac{1}{2}\int_{0}^{\infty}d\sigma_{1}\int_{0}^{\infty}d\sigma_{2}\frac{1}{\left(\sigma_{1}+1\right)^{3/2}}\frac{2}{\sigma_{2}^{3/2}}K_{0}\left(2\frac{\sqrt{a\sigma_{2}^{2}+b\sigma_{2}+c}}{\sqrt{\sigma_{2}}}\right) & =\\
\int d^{3}x_{2}\frac{1}{2}\frac{1}{2}\int_{0}^{\infty}d\sigma_{1}\frac{1}{\left(\sigma_{1}+1\right)^{3/2}}\frac{\pi e^{-2\sqrt{2\sqrt{a}\sqrt{c}+b}}}{\sqrt{c}}= & \Rightarrow\\
\int d^{3}x_{2}\frac{1}{2}\frac{1}{2}\int_{0}^{\infty}d\sigma_{1}\frac{1}{\left(\sigma_{1}+1\right)^{3/2}}\frac{\pi}{\sqrt{c}}e^{-2\left(\sqrt{a\left(\sigma_{1}+1\right)}+\sqrt{\frac{c}{\sigma_{1}+1}}\right)}
\end{array}\label{eq:result_1}
\end{equation}

The second term has the same general form, but with $\eta_{1}^{2}\leftrightarrow\eta_{12}^{2}$.
This means that we can write it as

\begin{equation}
\begin{array}{ccc}
\int d^{3}x_{2}\frac{1}{2}\frac{1}{2}\int_{0}^{\infty}d\sigma_{1}\int_{0}^{\infty}d\sigma_{2}\frac{2}{\left(\sigma_{1}+1\right)^{3/2}\sigma_{2}^{3/2}}K_{0}\left(\frac{2\sqrt{\frac{\eta_{13}^{2}\sigma_{2}^{2}x_{2}^{2}}{4\left(\sigma_{1}+1\right)}+\frac{1}{4}\left(\eta_{1}^{2}+\eta_{12}^{2}\sigma_{1}\right)x_{2}^{2}+\frac{\left(\eta_{1}^{2}+\eta_{12}^{2}\sigma_{1}+\eta_{13}^{2}\left(\sigma_{1}+1\right)\right)\sigma_{2}x_{2}^{2}}{4\left(\sigma_{1}+1\right)}}}{\sqrt{\sigma_{2}}}\right) & =\\
\int d^{3}x_{2}\frac{1}{2}\frac{1}{2}\int_{0}^{\infty}d\sigma_{1}\int_{0}^{\infty}d\sigma_{2}\frac{1}{\left(\sigma_{1}+1\right)^{3/2}}\frac{2}{\sigma_{2}^{3/2}}K_{0}\left(2\frac{\sqrt{a\sigma_{2}^{2}+\left(b-f\right)\sigma_{2}+\left(c-g\right)}}{\sqrt{\sigma_{2}}}\right) & =\\
\int d^{3}x_{2}\frac{1}{2}\frac{1}{2}\int_{0}^{\infty}d\sigma_{1}{\displaystyle \frac{1}{\left(\sigma_{1}+1\right)^{3/2}}\frac{\pi e^{-2\sqrt{2\sqrt{a}\sqrt{c-g}+\left(b-f\right)}}}{\sqrt{c-g}}}= & \Rightarrow\\
\int d^{3}x_{2}\frac{1}{2}\frac{1}{2}\int_{0}^{\infty}d\sigma_{1}\frac{1}{\left(\sigma_{1}+1\right)^{3/2}}\frac{\pi}{\sqrt{c}}e^{-2\left(\sqrt{a\left(\sigma_{1}+1\right)}+\sqrt{\frac{c-g}{\sigma_{1}+1}}\right)}
\end{array}\label{eq:result_2}
\end{equation}
where 

\begin{equation}
\left\{ \, g\to\frac{1}{4}x_{2}^{2}\left(\eta_{1}^{2}-\eta_{12}^{2}\right)\left(\sigma_{1}-1\right),f\to\frac{x_{2}^{2}\left(\eta_{1}^{2}-\eta_{12}^{2}\right)\left(\sigma_{1}-1\right)}{4\left(\sigma_{1}+1\right)}\right\} \label{eq:f_g}
\end{equation}

The third term has $x^{1/2}$ instead of $x^{3/2}$ in the denominator
so we have to use a different sequence of derivatives of (\ref{eq:PBM3 p. 349 No. 2.24.2.9})
to derive

\begin{align}
\int_{0}^{\infty}dx\,\frac{1}{x^{1/2}}K_{0}\left(2\sqrt{\frac{ax^{2}+bx+c}{x}}\right) & =\int_{0}^{\infty}dx\,\frac{1}{x^{1/2}}\sqrt{\pi}e^{-2\sqrt{\frac{ax^{2}+bx+c}{x}}}U\left(\frac{1}{2},1,4\sqrt{\frac{ax^{2}+bx+c}{x}}\right)\nonumber \\
 & =\int_{0}^{\infty}dx\,\frac{1}{2x^{1/2}}G_{0,2}^{2,0}\left(\frac{\text{ax}^{2}+\text{b x}+c}{x}|\begin{array}{c}
0,0\end{array}\right)\nonumber \\
 & =\frac{\sqrt{\pi}}{\sqrt{a}}G_{1,3}^{3,0}\left(b+2\sqrt{a}\sqrt{c}|\begin{array}{c}
-\frac{1}{2}\\
0,-\frac{1}{2},\frac{1}{2}
\end{array}\right)\quad,\label{eq:K0/x^1/2}\\
 & =\frac{e^{-2\sqrt{b+2\sqrt{a}\sqrt{c}}}\pi}{\sqrt{a}}\nonumber \\
 & \Rightarrow\frac{\pi}{\sqrt{a_{3}}}e^{-2\left(\sqrt{\frac{a_{3}}{\sigma_{1}+1}}+\sqrt{c_{3}\left(\sigma_{1}+1\right)}\right)}\nonumber 
\end{align}
where the reduction in the fourth line is from \cite{functions.wolfram.com/07.34.03.0955.01}
and the last line holds for the class of parameters akin to this particular
set: $\left\{ a_{3}\to\frac{1}{4}x_{2}^{2}\left(\eta_{1}^{2}\sigma_{1}+\eta_{12}^{2}\right),\, b_{3}\to\frac{x_{2}^{2}\left(\eta_{1}^{2}\sigma_{1}+\eta_{13}^{2}\left(\sigma_{1}+1\right)+\eta_{12}^{2}\right)}{4\left(\sigma_{1}+1\right)},\, c_{3}\to\frac{x_{2}^{2}\eta_{13}^{2}}{4\left(\sigma_{1}+1\right)}\right\} $.
We furthermore see that $\left\{ a_{3}\to c,\, c_{3}\to a\right\} $
in the present case, so the third term integrates to the same value
as the first:

\begin{equation}
\begin{array}{ccc}
\int d^{3}x_{2}\frac{1}{2}\frac{1}{2}\int_{0}^{\infty}d\sigma_{1}\int_{0}^{\infty}d\sigma_{2}\frac{2}{\left(\sigma_{1}+1\right){}^{3/2}\sqrt{\sigma_{2}}}K_{0}\left(\frac{2\sqrt{\frac{1}{4}\left(\sigma_{1}\eta_{1}^{2}+\eta_{12}^{2}\right)\sigma_{2}^{2}x_{2}^{2}+\frac{\left(\sigma_{1}\eta_{1}^{2}+\eta_{12}^{2}+\eta_{13}^{2}\left(\sigma_{1}+1\right)\right)\sigma_{2}x_{2}^{2}}{4\left(\sigma_{1}+1\right)}+\frac{\eta_{13}^{2}x_{2}^{2}}{4\left(\sigma_{1}+1\right)}}}{\sqrt{\sigma_{2}}}\right) & =\\
\int d^{3}x_{2}\frac{1}{2}\frac{1}{2}\int_{0}^{\infty}d\sigma_{1}\int_{0}^{\infty}d\sigma_{2}\frac{1}{\left(\sigma_{1}+1\right){}^{3/2}}\frac{2}{\sqrt{\sigma_{2}}}K_{0}\left(\frac{2\sqrt{c\sigma_{2}^{2}+b\sigma_{2}+a}}{\sqrt{\sigma_{2}}}\right) & =\\
\int d^{3}x_{2}\frac{1}{2}\frac{1}{2}\int_{0}^{\infty}d\sigma_{1}{\displaystyle \frac{1}{\left(\sigma_{1}+1\right){}^{3/2}}\frac{\pi e^{-2\sqrt{2\sqrt{a}\sqrt{c}+b}}}{\sqrt{c}}}= & \Rightarrow\\
\int d^{3}x_{2}\frac{1}{2}\frac{1}{2}\int_{0}^{\infty}d\sigma_{1}\frac{1}{\left(\sigma_{1}+1\right){}^{3/2}}\frac{\pi}{\sqrt{c}}e^{-2\left(\sqrt{a\left(\sigma_{1}+1\right)}+\sqrt{\frac{c}{\sigma_{1}+1}}\right)}
\end{array}\label{eq:result_3}
\end{equation}

The fourth term has the same general form as the third, but with $\eta_{1}^{2}\leftrightarrow\eta_{12}^{2}$.
This means that we can write it as

\begin{equation}
\begin{array}{ccc}
\int d^{3}x_{2}\frac{1}{2}\frac{1}{2}\int_{0}^{\infty}d\sigma_{1}\int_{0}^{\infty}d\sigma_{2}\frac{2}{\left(\sigma_{1}+1\right)^{3/2}\sigma_{2}^{1/2}}K_{0}\left(\frac{2\sqrt{\frac{1}{4}\left(\eta_{1}^{2}+\eta_{12}^{2}\sigma_{1}\right)\sigma_{2}^{2}x_{2}^{2}+\frac{\left(\eta_{1}^{2}+\eta_{12}^{2}\sigma_{1}+\eta_{13}^{2}\left(\sigma_{1}+1\right)\right)\sigma_{2}x_{2}^{2}}{4\left(\sigma_{1}+1\right)}+\frac{\eta_{13}^{2}x_{2}^{2}}{4\left(\sigma_{1}+1\right)}}}{\sqrt{\sigma_{2}}}\right) & =\\
\int d^{3}x_{2}\frac{1}{2}\frac{1}{2}\int_{0}^{\infty}d\sigma_{1}\int_{0}^{\infty}d\sigma_{2}\frac{1}{\left(\sigma_{1}+1\right)^{3/2}}\frac{2}{\sigma_{2}^{1/2}}K_{0}\left(2\frac{\sqrt{a\sigma_{2}^{2}+\left(b-f\right)\sigma_{2}+\left(c-g\right)}}{\sqrt{\sigma_{2}}}\right) & =\\
\int d^{3}x_{2}\frac{1}{2}\frac{1}{2}\int_{0}^{\infty}d\sigma_{1}{\displaystyle \frac{1}{\left(\sigma_{1}+1\right)^{3/2}}\frac{\pi e^{-2\sqrt{2\sqrt{a}\sqrt{c-g}+\left(b-f\right)}}}{\sqrt{c-g}}}= & \Rightarrow\\
\int d^{3}x_{2}\frac{1}{2}\frac{1}{2}\int_{0}^{\infty}d\sigma_{1}\frac{1}{\left(\sigma_{1}+1\right)^{3/2}}\frac{\pi}{\sqrt{c}}e^{-2\left(\sqrt{a\left(\sigma_{1}+1\right)}-\sqrt{\frac{c-g}{\sigma_{1}+1}}\right)}
\end{array}\label{eq:result_4}
\end{equation}
where $f$ and $g$ are as in (\ref{eq:f_g}).

Numerically integrating these final four analytical results, the final
lines of (\ref{eq:result_1}), (\ref{eq:result_2}), (\ref{eq:result_3}),
and (\ref{eq:result_4}), for the values of $\left\{ \eta_{1}\to0.3,\eta_{12}\to0.5,\eta_{13}\to0.9\right\} $ used
above, the four terms differ only in the final two decimal places

\begin{eqnarray}
S_{1}^{0.3,0,0.5,0,0.9,0}\left(0,0;0,0,0\right) & = & \left(29.37382253070293   \;+29.373822530702924+29.37382253070293 \right. \nonumber  \\
& +  &  \left. \hspace{0.1cm}  \, 29.373822530702924\right)=117.49529012281171\quad,\label{eq:S3_final_an}
\end{eqnarray}

\noindent
and the pentultimate lines in each case are identically $29.373822530702924$.
So our hope that the sum of the four integrals over the intervals
$\left[0,1\right]$ and $\left[1,\infty\right]$ could be analytically
found by using integration limits of $\left[0,\infty\right]$ as a
sort of Rosetta Stone bears fruit:

\begin{equation}
\begin{array}{ccc}
\int d^{3}x_{2}\int_{1}^{\infty}d\sigma_{1}\frac{1}{\left(\sigma_{1}+1\right)^{3/2}}\int_{1}^{\infty}d\sigma_{2}\left(\frac{2}{\sigma_{2}^{3/2}}\left[K_{0}\left(2\frac{\sqrt{a\sigma_{2}^{2}+b\sigma_{2}+c}}{\sqrt{\sigma_{2}}}\right)+K_{0}\left(2\frac{\sqrt{a\sigma_{2}^{2}+\left(b-f\right)\sigma_{2}+\left(c-g\right)}}{\sqrt{\sigma_{2}}}\right)\right]\right. & +\\
\left.\frac{2}{\sqrt{\sigma_{2}}}\left[K_{0}\left(\frac{2\sqrt{c\sigma_{2}^{2}+b\sigma_{2}+a}}{\sqrt{\sigma_{2}}}\right)+K_{0}\left(2\frac{\sqrt{a\sigma_{2}^{2}+\left(b-f\right)\sigma_{2}+\left(c-g\right)}}{\sqrt{\sigma_{2}}}\right)\right]\right) & =\\
31.4147+22.115+38.9735+24.9916=117.495 & =\\
\int d^{3}x_{2}\int_{1}^{\infty}d\sigma_{1}\frac{1}{\left(\sigma_{1}+1\right)^{3/2}}\left(\frac{\pi}{\sqrt{c}}e^{-2\sqrt{2\sqrt{a}\sqrt{c}+b}}+\frac{\pi}{\sqrt{c-g}}e^{-2\sqrt{2\sqrt{a}\sqrt{c-g}+\left(b-f\right)}}\right. & +\\
\left.\frac{\pi}{\sqrt{c}}e^{-2\sqrt{2\sqrt{a}\sqrt{c}+b}}+\frac{\pi}{\sqrt{c-g}}e^{-2\sqrt{2\sqrt{a}\sqrt{c-g}+\left(b-f\right)}}\right) & \Rightarrow & \quad.\\
\int d^{3}x_{2}\int_{1}^{\infty}d\sigma_{1}\frac{1}{\left(\sigma_{1}+1\right)^{3/2}}\left(\frac{\pi}{\sqrt{c}}e^{-2\left(\sqrt{a\left(\sigma_{1}+1\right)}-\sqrt{\frac{c}{\sigma_{1}+1}}\right)}+\frac{\pi}{\sqrt{c}}e^{-2\left(\sqrt{a\left(\sigma_{1}+1\right)}-\sqrt{\frac{c-g}{\sigma_{1}+1}}\right)}\right. & +\\
\left.\frac{\pi}{\sqrt{c}}e^{-2\left(\sqrt{a\left(\sigma_{1}+1\right)}+\sqrt{\frac{c}{\sigma_{1}+1}}\right)}+\frac{\pi}{\sqrt{c}}e^{-2\left(\sqrt{a\left(\sigma_{1}+1\right)}-\sqrt{\frac{c-g}{\sigma_{1}+1}}\right)}\right) & =\\
35.1943+23.5533+35.1943+23.5533=117.495
\end{array}\label{eq:4ints_an}
\end{equation}

We note that the analytical  $\sigma_{2}$  results (in the last line above after numerically integrating over $ d^{3}x_{2}d\sigma_{1}$) are identical for the first and
third terms, and for the second and fourth. One can likewise see that
the average of the first and third terms resulting from the numerical integral that includes $\sigma_{2}$ (in the third line above) 
gives the first analytical term, and the average of the
second and fourth terms resulting from the numerical integral that includes $d \sigma_{2}$  gives the second
analytical term. So our results extend significantly beyond
our hope that the sum of the four integrals over the interval $\left[1,\infty\right]$
would have an analytical result. If fact, an analytical result results
from the sum of only two terms. The same holds when we replace the
integral limits with $\left[0,1\right]$. Written with this new understanding,
the two pairs of integrals give

\begin{equation}
\begin{array}{ccc}
\int d^{3}x_{2}\int_{0}^{1}d\sigma_{1}\frac{1}{\left(\sigma_{1}+1\right)^{3/2}}\int_{0}^{1}d\sigma_{2}\left(\frac{2}{\sigma_{2}^{3/2}}K_{0}\left(2\frac{\sqrt{a\sigma_{2}^{2}+b\sigma_{2}+c}}{\sqrt{\sigma_{2}}}\right)+\frac{2}{\sqrt{\sigma_{2}}}K_{0}\left(\frac{2\sqrt{c\sigma_{2}^{2}+b\sigma_{2}+a}}{\sqrt{\sigma_{2}}}\right)\right) & =\\
\int d^{3}x_{2}\int_{0}^{1}d\sigma_{1}\frac{1}{\left(\sigma_{1}+1\right)^{3/2}}\frac{\pi}{\sqrt{c}}e^{-2\sqrt{2\sqrt{a}\sqrt{c}+b}} & \Rightarrow \\
\int d^{3}x_{2}\int_{0}^{1}d\sigma_{1}\frac{1}{\left(\sigma_{1}+1\right)^{3/2}}2\frac{\pi}{\sqrt{c}}e^{-2\left(\sqrt{a\left(\sigma_{1}+1\right)}-\sqrt{\frac{c}{\sigma_{1}+1}}\right)} & =\\
2\times23.5533\\
\\
\int d^{3}x_{2}\int_{0}^{1}d\sigma_{1}\frac{1}{\left(\sigma_{1}+1\right)^{3/2}}\int_{0}^{1}d\sigma_{2}\left(\frac{2}{\sigma_{2}^{3/2}}K_{0}\left(2\frac{\sqrt{a\sigma_{2}^{2}+\left(b-f\right)\sigma_{2}+\left(c-g\right)}}{\sqrt{\sigma_{2}}}\right)+\frac{2}{\sqrt{\sigma_{2}}}K_{0}\left(2\frac{\sqrt{a\sigma_{2}^{2}+\left(b-f\right)\sigma_{2}+\left(c-g\right)}}{\sqrt{\sigma_{2}}}\right)\right) & =\\
\int d^{3}x_{2}\int_{0}^{1}d\sigma_{1}\frac{1}{\left(\sigma_{1}+1\right)^{3/2}}2\frac{\pi}{\sqrt{c-g}}e^{-2\sqrt{2\sqrt{a}\sqrt{c-g}+\left(b-f\right)}} & \Rightarrow \\
\int d^{3}x_{2}\int_{0}^{1}d\sigma_{1}\frac{1}{\left(\sigma_{1}+1\right)^{3/2}}2\frac{\pi}{\sqrt{c}}e^{-2\left(\sqrt{a\left(\sigma_{1}+1\right)}-\sqrt{\frac{c-g}{\sigma_{1}+1}}\right)} & =\\
2\times35.1943 & =
\end{array}\label{eq:4ints[0,1]_an}
\end{equation}
whose sum is again $117.495$. 

We noted earlier that one can have a different $\sigma_{1}$ integration
interval than for $\sigma_{2}$. If we substitute $\left[0,\infty\right]$
for the $\sigma_{1}$ integration interval, and multiply by $\frac{1}{2},$
while keeping $\left[0,1\right]$ for the $\sigma_{2}$ integral,
each of the above two integrals yield a value of $2\times29.3738$
whose sum is again $117.495$. These are averages of the first and
third terms in set of results from the individual terms, \{39.0904,
43.4151, 19.6571, 15.3325\}. The same holds for we substitute $\left[0,\infty\right]$
for the $\sigma_{1}$ integration interval, and multiply by $\frac{1}{2},$
while setting $\left[1,\infty\right]$ for the $\sigma_{2}$ integral:
though the first and last pairs swap values in the individual terms. 

These are integral results extracted from the particular case under
study. Of course the final form will not hold for generic values of
the parameters, but do the penultimate ones hold? Indeed, they do:

\begin{equation}
\begin{array}{ccc}
\int_{0}^{1}\, dx\left(\frac{2}{x^{3/2}}K_{0}\left(\frac{2\sqrt{ax^{2}+bx+c}}{\sqrt{x}}\right)+\frac{2}{\sqrt{x}}K_{0}\left(\frac{2\sqrt{cx^{2}+bx+a}}{\sqrt{x}}\right)\right) & =\\
\int_{0}^{1}\, dx\left(\frac{2\sqrt{\pi}}{x^{3/2}}e^{-\frac{2\sqrt{ax^{2}+bx+c}}{\sqrt{x}}}U\left(\frac{1}{2},1,\frac{4\sqrt{ax^{2}+bx+c}}{\sqrt{x}}\right)+\frac{2\sqrt{\pi}}{\sqrt{x}}e^{-\frac{2\sqrt{a+bx+cx^{2}}}{\sqrt{x}}}U\left(\frac{1}{2},1,\frac{4\sqrt{cx^{2}+bx+a}}{\sqrt{x}}\right)\right) & =\\
\int_{0}^{1}\, dx\left(\frac{1}{x^{3/2}}G_{0,2}^{2,0}\left(\frac{ax^{2}+bx+c}{x}|\begin{array}{c}
0,0\end{array}\right)+\frac{1}{\sqrt{x}}G_{0,2}^{2,0}\left(\frac{cx^{2}+bx+a}{x}|\begin{array}{c}
0,0\end{array}\right)\right) & =\\
\frac{\pi}{\sqrt{c}}e^{-2\sqrt{2\sqrt{a}\sqrt{c}+b}} & =\\
0.738215\\
\left[\{a\to0.21,\, b\to0.31,\, c\to0.41\}\right]
\end{array}\label{eq:2ints_an_generic}
\end{equation}
The same values result when the integration is over the interval $\left[1,\infty\right]$,
as they do if we change the integration limits to $\left[0,\infty\right]$
and multiply by $\frac{1}{2}$. We note that the second pair of integrals
in (\ref{eq:4ints[0,1]_an}), when cast into generic terms, replicate
the above integral by simply renaming $b-f$ as $b$ and $c-g$ as
$g$, so we have derived one new integral relation and not two distinct ones.

While each if of the two terms in the above integral does not map
onto a known analytical result, having the sum map onto a known analytical
result is a significant step in that direction.

\section*{Conclusion}

We have crafted a quartet of integral representations of products
of Slater orbitals over the interval $[0,1]$ that may be useful in
the reduction of multidimensional transition amplitudes of quantum
theory. For three of these representations, the general form was found
for a product of any number of Slater orbitals, whose derivatives
in the atomic realm represent hydrogenic and Hylleraas wave functions,
as well as those composed of explicitly correlated exponentials of
the kind introduced by Thakkar and Smith.\cite{Thakkar and Smith77}
These results are also useful in nuclear transition amplitudes, and
may also find application in solid-state physics, plasma physics,
negative ion physics, and problems involving a hypothesized non-zero-mass
photon.

These three integral representations have the advantage over the Gaussian
transform of requiring the introduction of one fewer integral to be
subsequently reduced. They also require many fewer than the $\left(3\left(M-1\right)+M-1\right)$
integral dimensions that the Fourier transform introduces for a product
of M Slater orbitals. Direct integration of products of Slater orbitals
containing angular functions centered on different points, usually
bears fruit in only the simplest problems. The fourth conventional
approach to these problems is to represent M Slater orbitals as M
addition theorems, that is M infinite sums over Spherical Harmonics
containing the angular dependencies. Orthogonality allows one to remove
a few of these infinite sums in the process of integrating some of
the original integrals. For large M, this approach rapidly bogs down.

Each of the four extant approaches to such problems  run into
difficulties at some point as M increases. These three new integral
representations (over the interval $[0,1]$) for M Slater orbitals
likewise have some positives and some negatives. We have found that
only two of these three integral representations over the interval
$[0,1]$ for M Slater orbitals is numerically stable for $M>3$, though
the other might be better for analytical work in some cases. The prior
version having integrals running over the interval $[0,\infty]$ was
also numerically stable, but one does have the inconvenience
of needing to test for a sufficiently large upper integration limit.

For the simplest problems, the three integral representations (over
the interval $[0,1]$) for M Slater orbitals of the present paper
provide solutions in a much more rapid fashion than do the four extant
approaches, and even surpass the integral representation over the
interval$[0,\infty]$ of the prior paper in allowing the moderately
hard problem of the integral over three Slater orbitals -- after all
coordinate integrals have been done -- to be reduced to analytical
form via tabled integrals over the interval $[0,1]$. On the other
hand, the integral representation over the interval $[0,\infty]$ of
the prior paper lacked any tabled result in the final step of this problem
and had to rely on the computer algebra and calculus program Mathematica 7 to do this integral. This, however, belies the general paucity of
tabled integrals over the interval $[0,1]$ relative to those over
the interval $[0,\infty]$.

The fact that the integration variables reside within a square root
as the argument of a Macdonald function, shared by the prior work,
will lead to difficulties in some complicated problems since only
one such integral (transformed into a Meijer G-function) was known
prior to the previous paper. Unlike that paper, the three integral
representations for M Slater orbitals of the present work have the
added difficulty of having no known integrals of this sort upon which
to build. 

It was for this reason that we introduced a fourth integral representation
that is not easily generalizable to large M, but one hoped it would provide
a bridge for finding the requisite integrals in the above problems. This final integral
representation allowed us to derive the analytical result for an integral
of a sum of two Meijer G functions \textmd{\normalsize{} $f\left(x\right)G_{0,2}^{2,0}\left(\frac{\text{ax}^{2}+\text{b x}+c}{x}|\protect\begin{array}{c}
0,0\protect\end{array}\right)$} over the interval $[0,1]$, via
a bridge from the version of this integral representation that is over the
interval $[0,\infty]$. This is only half-way to the desired result,
but is a promising step, and provides an integral that researchers
in fields far afield from atomic theory may find useful.

\section*{Acknowledgements}

As it became clear that I needed to make at least some progress on finding analytical results for integrals
of Meijer G-functions $f\left(x\right)G_{0,2}^{2,0}\left(\frac{\text{ax}^{2}+\text{b x}+c}{x}|\begin{array}{c}
0,0\end{array}\right)$ over the interval  $\left[0,1\right]$ for this  research project to come to a sense of completion, I chanced upon an integral in Gr\"{o}bner und Hofreiter\cite{GH p. 176 No. 421.8} that would serve as an integral representation for a product of denominators.  I found that I could use this as the basis for an integral representation for a product of several Slater orbitals that had the property of bridging from known integrals of Meijer G-functions (with such arguments)  over the interval  $\left[0,\infty\right]$ to heretofore untabled pairs of such integrals over the
 the interval  $\left[0,1\right]$ (and over $\left[1,\infty\right]$); a sort of mathematical Rosetta Stone.

I have long had the practice of filling the 10 minutes prior to when my Astronomy class starts with videos of music featuring women instrumentalists, just as I make it my practice to bring video clips of experts in the field who happen to be women into the class content. The reader may or may not be aware that the historical predicament
of women in STEM fields has significant parallels to the historical
predicament of women in music, particularly when it comes to women
instrumentalists. Indeed, women came to be auditioned into orchestras
in significant numbers only after blind auditions were introduced
in the 1970s. Neither did one see women instrumentalists playing with
Miles Davis, say, or The Rolling Stones. Fortunately, both pop music
and STEM fields are beginning to shift in this regard. It is my hope that as the younger generation comes to see women in both roles as ``normal," they will help accelerate this shift.

So my students hear Lari Basilio shredding on electric guitar,\cite{Basilio} Sonah Jobarteh on the kora,\cite{Jobarteh} and  Sophie Alloway on drums with Ida Hollis on electric bass, among many others.  They likewise  learn about the process of looking for life on Mars from Dr. Moogega Cooper,\cite{Moogega} and about the sound of Back Holes colliding from Dr. Janna Levin.\cite{Levin} 

I share all of this detail so that it will be clear why the background soundtrack to my research into the material that comprises Sections 7 and 8 -- and the idea that one sort of integral could act as a bridge or mathematical Rosetta Stone to craft others -- was guitarist, composer, and singer Sister
Rosetta Tharpe.\cite{Tharpe}  She was  inducted by The Rock and Roll Hall
of Fame as  \textquotedblleft{}the Godmother of Rock \&
Roll,\textquotedblright{}\cite{rockhall_re_tharpe} though her 
influences on gospel, country, and R\&B were also vast. 

I am a
jazz drummer who has been immersed in learning these other four musical
styles over the past two decades, and, thus, Sister Rosetta Tharpe has been
key to not only to deepening and generalizing my musical patterning, but also
to the joy I experience in the process. 
On a weekly basis I am immersed in the creative expression of the
musicians I jam with, and the rhythmical patterns they manifest in
their music evoke a resonant rhythmical response in my drumming. Sometimes
this response is delayed by months, because my skill-level needs to
grow to accommodate it. And I sense, but cannot prove, that this response also manifests in my work as  a theoretical physicist who relies heavily on
pattern recognition for insights that culminate in my math-based results,
such as the generalization to eq. (\ref{eq:MtransCompact_sig}) from
the sequence from (\ref{eq:3-sig}) to (\ref{eq:5-sig}). 

It is with all this in mind that I dedicate this paper, and in particuar the integral representations of Sections 7 and 8,  to Sister
Rosetta Tharpe.

\end{document}